\documentclass[a4paper]{amsart}
\usepackage{amssymb}
\usepackage[hidelinks]{hyperref}
\usepackage{tikz}
\usepackage{graphicx}
\usepackage{todonotes}
\usepackage[font=small]{caption}
\usetikzlibrary{decorations.markings}

\newcommand{\Z}{\mathbb{Z}}
\newcommand{\R}{\mathbb{R}}
\newcommand{\C}{\mathbb{C}}
\newcommand{\T}{\mathbb{T}}
\newcommand{\ii}{{\rm i}}
\newcommand{\cent}[1]{{#1}^\circ}

\newcommand{\surface}{\mathcal}
\newcommand{\sM}{\surface{M}}
\newcommand{\sS}{\surface{S}}
\newcommand{\bx}{\mathbf{x}}
\newcommand{\ff}{\mathbf{f}}
\newcommand{\bT}{\mathbf{T}}
\newcommand{\tm}{\text{max}}

\newcommand{\re}{\operatorname{Re}}
\newcommand{\im}{\operatorname{Im}}
\newcommand{\Res}[2]{\operatorname{Res}\left(#1,#2\right)}

\theoremstyle{plain} 
\newtheorem{theorem}{Theorem}
\newtheorem{lemma}[theorem]{Lemma}

\newtheorem{proposition}[theorem]{Proposition}
\theoremstyle{definition} 

\theoremstyle{remark}
\newtheorem{remark}{Remark}

\begin{document}

\title[Gluing doubly Scherk]{Gluing Doubly Periodic Scherk Surfaces\\ into Minimal Surfaces}

\author{Hao Chen}
\email{chenhao5@shanghaitech.edu.cn}
\author{Yunhua Wu}
\email{wuyh12023@shanghaitech.edu.cn}
\address{ShanghaiTech University, Shanghai, China}


\date{\today}

\begin{abstract}
	We construct minimal surfaces by stacking doubly periodic Scherk surfaces one
	above another and gluing them along their ends.  It is previously known that
	the Karcher--Meeks--Rosenberg (KMR) doubly periodic minimal surfaces and Meeks'
	family of triply periodic minimal surfaces can both be obtained by gluing two
	Scherk surfaces.  There have been hope and failed attempts to glue more Scherk
	surfaces.  But our analysis shows that: Except for the special case where the
	doubly periodic Scherk surfaces all have triangular horizontal lattice, a glue
	construction can only produce the trivial Scherk surface itself, the KMR
	examples, or Meeks' surfaces.
\end{abstract}

\maketitle

\section{Introduction}

The goal of this paper is to glue finitely many doubly periodic Scherk surfaces
into minimal surfaces.

\emph{Doubly periodic Scherk surfaces} (or simply \emph{Scherk surfaces}
hereafter) form a $1$-parameter family of doubly periodic minimal surfaces
(DPMSs).  If we quotient out the periods, a Scherk surface $\sS$ is a minimal
surface of genus $0$ in $\T^2 \times \R$ with four vertical Scherk ends
(asymptotic to vertical planes).  Two of the ends extend upwards and both are
parallel to $\bT_1 = (\cos\theta_1, \sin\theta_1)$, and the other two extend
downwards and both are parallel to $\bT_2 = (\cos\theta_2, \sin\theta_2)$.
Here, $\T^2$ is a flat torus whose fundamental parallelogram is a rhombus
spanned by $2\pi\bT_1$ and $2\pi\bT_2$.  $\sS$ is a minimal graph over the
domain
\[
  \{a \bT_1 + b \bT_2 \mid (a,b) \in (0, \pi)^2 \cup (\pi, 2\pi)^2\}.
\]
So it looks like a rhombic checkerboard when seen from above and $\sS$ is only
defined over the ``black'' rhombi; see Figure~\ref{fig:topview}.  From far away,
$\sS$ looks like vertical planes parallel to $\bT_1$ suddenly transition into
vertical planes parallel to $\bT_2$.  The Gaussian curvature accumulates near a
horizontal plane where the transition occurs.

\begin{figure}[hbt]
	\includegraphics[height=0.4\textwidth]{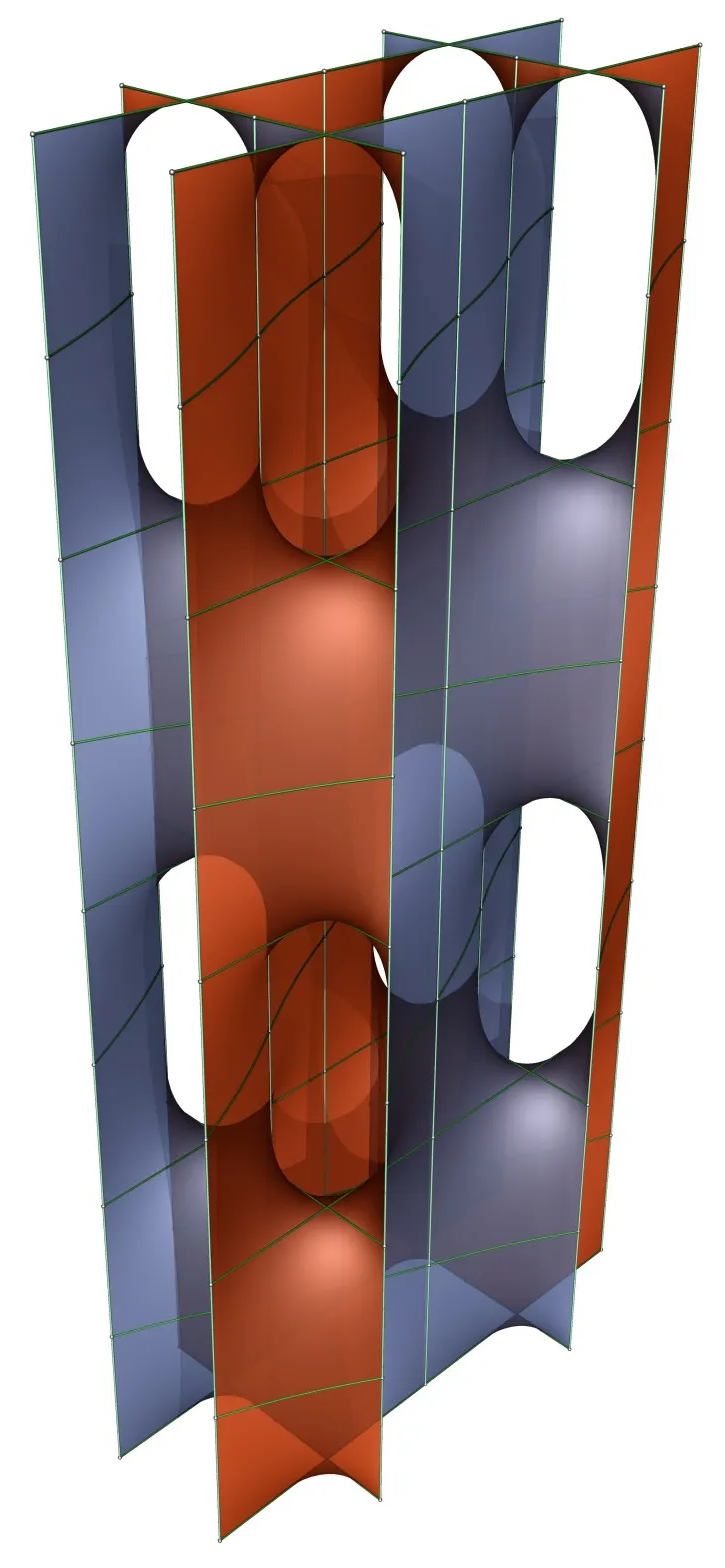}
	\includegraphics[height=0.4\textwidth]{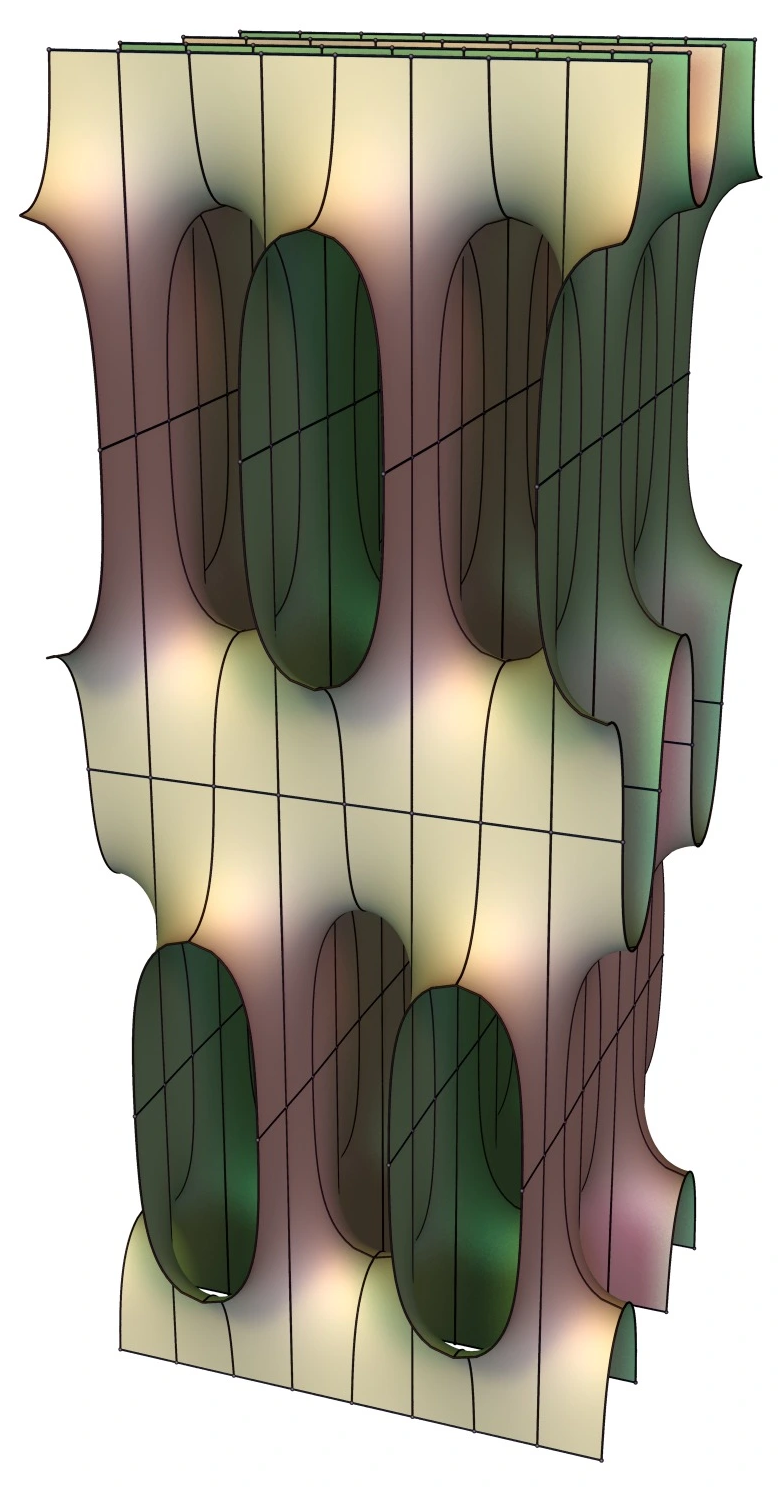}
	\includegraphics[height=0.4\textwidth]{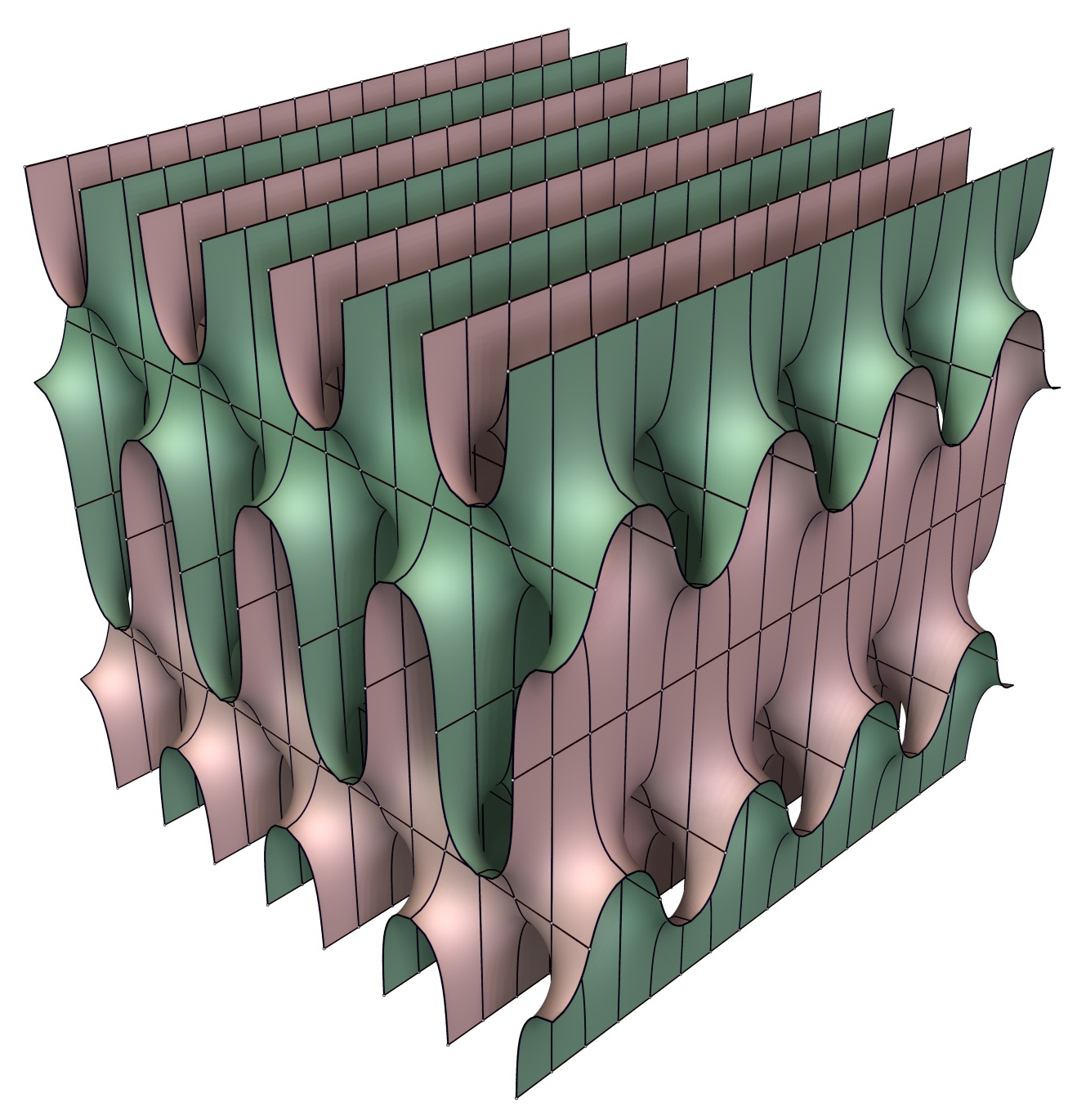}
	\caption{
		Triply periodic CLP (left), oP (middle), and oD (right) surfaces near the
		doubly periodic Scherk limit.
	}
	\label{fig:meeks}
\end{figure}

Karcher--Meeks--Rosenberg (KMR) surfaces are DPMSs of genus $1$ with four ends
that form a three-parameter family.  They admit limits that look like two Scherk
surfaces glued along their ends.  Meeks' family of triply periodic minimal
surfaces (TPMSs) of genus 3 also admit limits that look like two Scherk surfaces
glued together.  For instance, Figure~\ref{fig:meeks} shows Scherk limits of
Schwarz' CLP, oP, and oD families.  These examples stimulate hope and (failed)
attempts to glue more Scherk surfaces into minimal surfaces of higher genus.

We will show that, except for the special case of triangular horizontal lattice,
one can not construct a periodic minimal surface by gluing more than two but
finitely many Scherk surfaces.  More precisely
\begin{theorem}\label{thm:main}
  Let $\sM_\varepsilon$, $\varepsilon>0$, be a one-parameter family of
  minimal surfaces in $\T_\varepsilon^2 \times \R$ such that, as $\varepsilon
  \to 0$ and for any $H \in \R$, $\sM_\varepsilon - (0, 0, H /
  \varepsilon^2)$ converges either to a doubly periodic Scherk surface, or to
  a pair of vertical planes.  Then
  \begin{itemize}
    \item either the horizontal flat torus $\T_\varepsilon^2$ converges to a
      triangular lattice as $\varepsilon \to 0$, 
    \item or otherwise, for $\varepsilon$ sufficiently small,
      $\sM_\varepsilon$ is either a doubly periodic Scherk surface, a KMR
      example, or a member of the Meeks family.
  \end{itemize}
\end{theorem}

We must stress that, by a ``doubly periodic Scherk surface in $\T^2 \times
\R$'', we only mean the DPMS of genus $0$ with four Scherk ends. Otherwise, one
could always choose a larger horizontal flat torus to achieve a higher genus.
This will bring more flexibility to the gluing construction.  This is however
not considered in the current paper.

Then the triangular lattice turns out to be the only horizontal lattice that
could produce higher-genus examples.  In fact, a recent numeric study has
constructed a TPMS of genus $4$ with a limit that looks like three Scherk
surfaces glued together, and the horizontal lattice is indeed
triangular~\cite{markande2018}. A rigorous proof for this example is under
preparation.  Even more complicated examples is possible with the triangular
horizontal lattice; see Remark~\ref{rem:triangle}. However, this special case
is not the focus of the current paper.  We will mainly talk about the
"otherwise" part of Main Theorem~\ref{thm:main}.

This paper is one in a series that aims to completely describe the ideal
boundary of the moduli space of the TPMSs of genus 3.  Previously, TPMSs of
genus 3 have been constructed by gluing catenoids, helicoids, and simply
periodic Scherk surfaces~\cite{traizet2008, chen2022, chen2024}.

\subsection{Glue Scherk surfaces into TPMSs}

We first try to construct TPMSs that look like $n$ Scherk surfaces, say $\sS_1,
\sS_2, \cdots, \sS_n$, stacked one above another and \emph{periodically} glued
along their ends.  More precisely, for $1 \le k \le n$, the upward ends of
$\sS_k$ are glued with the downward ends of $\sS_{k+1}$, where the subscript is
taken modulo $n$ (same afterwards).  It is then necessary that $n=2m$ is even
and (the asymptotic planes of) the upward ends of $\sS_k$ and the downward ends
of $\sS_{k+1}$ are parallel.

In the current paper, we only consider the simple case where the ends take only
two directions, depending on the parity of $k$.  For this purpose, we define the
parity function
\[
	\varsigma(k) = \begin{cases}
    1, & k \text{ odd},\\
    2, & k \text{ even},
	\end{cases}
\]
and we assume that the upward (resp.\ downward) ends of $\sS_k$ are parallel to
$\bT_{\varsigma(k)}$ (resp.\ $\bT_{\varsigma(k-1)}$).  This assumption seems to
be an oversimplification, but is actually generally true.  We will prove that
\begin{lemma}\label{lem:triangle}
  Under the assumptions of Theorem~\ref{thm:main},
  \begin{itemize}
    \item either the horizontal flat torus $\T_\varepsilon^2$ converges to
      a triangular lattice as $\varepsilon \to 0$, and the limit vertical
      planes have three possible directions,
    \item or otherwise, the vertical planes alternate between two possible
      directions.
  \end{itemize}
\end{lemma}

The surfaces we plan to construct form families $\sM_\varepsilon$ parameterized
by a real parameter $\varepsilon$.  At the beginning, $\varepsilon = 0$ and the
$\sS_k$'s are infinite distance apart.  Between $\sS_k$ and $\sS_{k+1}$, the
surface $\sM_0$ appears as vertical planes parallel to $\bT_{\varsigma(k)}$.
After a scaling by $\varepsilon^2$, the distances between the $\sS_k$'s become
finite, but the vertical planes become infinitesimally close. The Gaussian
curvature explodes at horizontal planes where the directions of the vertical
planes suddenly change.  Let $h_k$ denote the heights of these horizontal
planes, which can be seen as the height of $\sS_k$.  Then $\varepsilon^2
\sM_\varepsilon$ looks like foliations of $\T^2 \times (h_k, h_{k+1})$ by
vertical planes.  We define $\ell_k = h_{k+1} - h_k$, $1 \le k \le n$, as the
scaled distance between the $\sS_k$ and $\sS_{k+1}$, and write $\ell =
(\ell_k)_{1 \le k \le n}$ as a finite sequence of length $n$.  Up to a
reparameterization and a rescaling, we may assume that $\sum_k \ell_k = 1$.

\begin{figure}[hbt]
	\includegraphics[width=0.45\textwidth]{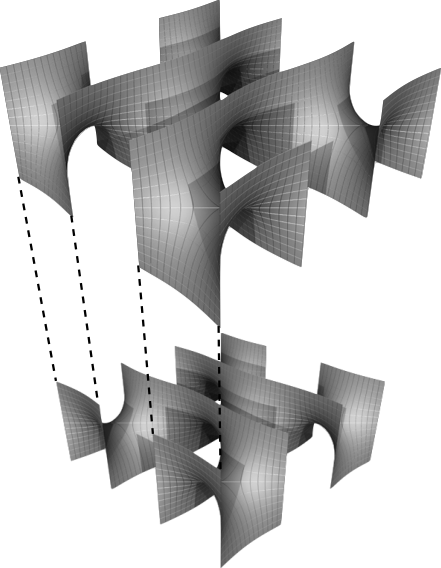}
	\includegraphics[width=0.45\textwidth]{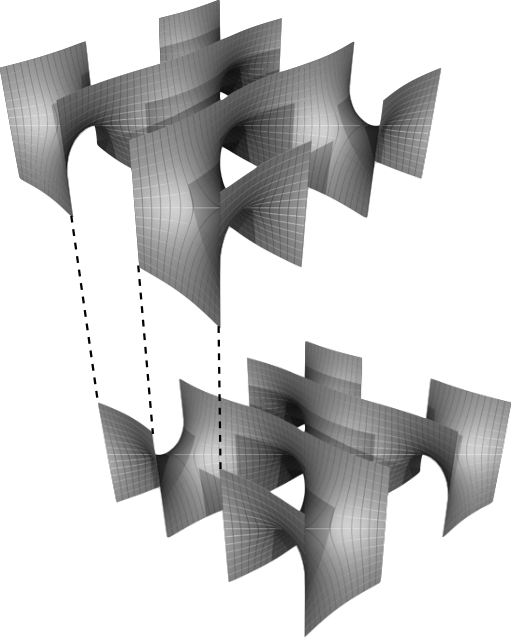}
	\caption{
    Two doubly periodic Scherk surfaces glued at phase $0$ (left) and phase
    $\pi$ (right).
	}
	\label{fig:phase}
\end{figure}

Moreover, it is possible that $\sS_{k+1}$ ``slides'' with respect to $\sS_k$
along the glued end.  The amount of the sliding, termed ``phase difference''
and denoted $\psi_k$, is defined as follows:  When the domain of $\sS_{k+1}$ as
a minimal graph coincides with the domain of $\sS_k$, the phase difference of
$\sS_{k+1}$ with respect to $\sS_k$ is $\psi_k = 0$, and we say $\sS_{k+1}$ is
\emph{in phase} with $\sS_k$; see Figure~\ref{fig:phase} (left).  In general,
$\sS_{k+1}$ is translated from the in-phase position by a vector
$\psi_k\bT_{\varsigma(k)}$, where $\psi_k \in \R/2\pi\Z$ is the phase
difference.  Note that, when $\psi_k = \pi$, the closures of the domains of
$\sS_k$ and $\sS_{k+1}$ cover the plane; in this case, we say that $\sS_{k+1}$
is \emph{in opposite phase} with $\sS_k$; see Figure~\ref{fig:phase} (right).
The sequence $\psi = (\psi_k)_{1 \le k \le n}$ is again a finite sequence of
length $n$.

The pair $(\ell, \psi)$ is termed a \emph{configuration}.  They give a complete  
description of the gluing pattern.

As $\varepsilon$ increases, the idea of the gluing is to desingularize the
initial surface.  There are two ways to look at the desingularization.  Either
the $\sS_k$'s are pulled a little bit closer towards each other, hence the
distances between them are reduced.  Alternatively, after a scaling by
$\varepsilon^2$, the vertical planes are pushed a little bit away from each
other, hence the horizontal periods increase.

\medskip

Our construction allows the horizontal flat torus to differ from rhombic after
desingularization.   Hence the horizontal periods, initially $2\pi\bT_1$ and
$2\pi\bT_2$ for $\varepsilon=0$, may not be unit vectors for $\varepsilon>0$.
Define
\[
  \ell_\tm = \max_k \ell_k,\qquad
  \tau_\tm = \exp(-\ell_\tm/\varepsilon^2).
\]
We assume that the horizontal periods are of the form
\begin{equation}\label{eq:horper}
	2\pi\bT_1 (1 + \Lambda_1 \tau_\tm), \qquad
	2\pi\bT_2 (1 + \Lambda_2 \tau_\tm),
\end{equation}
where $\Lambda_1$ and $\Lambda_2$ are real numbers.  See the proof of
Proposition~\ref{prop:alpharho} for the motivation of this assumption.  Note
that, as $\tau \to 0$, the aspect ratio is
\[
  \frac{1 + \Lambda_1 \tau_\tm}{1 + \Lambda_2 \tau_\tm} = 1 + (\Lambda_1 - \Lambda_2) \tau_\tm + o(\tau_\tm).
\]
So the aspect ratio is determined to the first order by $\Lambda_1 -
\Lambda_2$.  Adding the same constant to both $\Lambda$ corresponds to a
scaling of the horizontal lattice to the first order.

Define
\[
  \Psi_1 = \sum_{j=1}^{n/2} \psi_{2j-1},\qquad
  \Psi_2 = \sum_{j=1}^{n/2} \psi_{2j},\qquad
  \bT_3 = \Psi_1 \bT_1 + \Psi_2 \bT_2,
\]
Then $\bT_3 + (\Psi_1\Lambda_1\bT_1 + \Psi_2\Lambda_2\bT_2)\tau_\tm$ is the
horizontal component of the "vertical" period, while $1/\varepsilon^2$ is the
horizontal component.

To summarise: Given a periodic configuration $(\ell, \psi)$ of length $n = 2m$.
We expect, for sufficiently small $\varepsilon > 0$, a family $\sM_\varepsilon$
of minimal surfaces of genus $n+1$, in the flat torus
\begin{multline}
	\T_\varepsilon^3 = \R^3 / 2\pi\langle
	(\bT_1 (1 + \Lambda_1 \tau_\tm), 0), (\bT_2 (1 + \Lambda_2 \tau_\tm), 0),\\
	(\bT_3 + (\Psi_1 \Lambda_1 \bT_1 + \Psi_2 \Lambda_2 \bT_2) \tau_\tm, 1/\varepsilon^2)
	\rangle.
\end{multline}
They lift to TPMSs $\tilde\sM_\varepsilon$ in $\R^3$ such that
\begin{itemize}
	\item For each $1 \le k \le n$, there exists functions $H_k(\varepsilon)$ that
		diverges to infinity while $\sM_\varepsilon - (0,0,H_k(\varepsilon))$ converges
		to a doubly periodic Scherk surface $\sS_k$ as $\varepsilon \to 0$.  Moreover,
		$\varepsilon^{-2}(H_{k+1} - H_k) \to \ell_k$ as $\varepsilon \to 0$.

	\item As $\varepsilon\to 0$, $\tilde\sM_\varepsilon$ converges, after a
		scaling by $\varepsilon^2$, to foliations by vertical planes of the direct
		product of $\R^2$ with segments of length $\ell_k$.

	\item As $\varepsilon \to 0$, the phase difference between $\sS_{k+1}$ and
		$\sS_k$ converges to $\psi_k$.
\end{itemize}

In general, for such a gluing construction to succeed, a balance condition is
necessary, and a non-degenerate condition is needed for the use of the Implicit
Function Theorem.  In the current construction, the balance condition requires
that $\ell_k$ and $\psi_k$ only depend on the parity of $k$.  That is,
\[
  (\ell_k, \psi_k) = (\ell_{\varsigma(k)}, \psi_{\varsigma(k)}).
\]
We will also require that
\[
  K_k = \begin{cases}
    4\cos(\psi_k) + \Lambda_{\varsigma(k)} \exp(\zeta_k),& \ell_k = \ell_\tm,\\
    4\cos(\psi_k),& \ell_k < \ell_\tm,
  \end{cases}
\]
is nonzero for all $1 \le k \le n$, where $\zeta_k = \partial \ell_k /
\partial(\varepsilon^2)$.  This condition is not here to allow the use of the
Implicit Function Theorem, but it does allow us to proceed with the construction.

Our main result about TPMSs is the following
\begin{theorem}\label{thm:Meeks}
	For the family $\sM_\varepsilon$ of triply periodic minimal surfaces
	described above to exist, it is necessary that
	\[
    (\ell_k, \psi_k) = (\ell_{\varsigma(k)}, \psi_{\varsigma(k)}).
	\]
	If this is indeed the case and $K_k \ne 0$ for all $1 \le k \le n$, then
	$\ell$ is necessarily constant, and $\sM_\varepsilon$ does exist for
	sufficiently small $\varepsilon$ and coincides with the Meeks' family of
	triply periodic minimal surfaces of genus $3$.
\end{theorem}

The Meeks' family form a 5-dimensional real manifold, which is compatible to
the 5-dimensional deformation space of the lattice.  However, near the Scherk
limit, the manifold is not parameterized by the shapes of the lattice.  We will
see in Remark~\ref{rm:nosolution} that $\zeta_{1,2}$ replace $\Lambda_{1,2}$ as
parameters of the Meeks' manifold.

\subsection{Glue Scherk surfaces into DPMSs}

Now we try to stack only a finite sequence of Scherk surfaces and construct
DPMSs, leaving four free Scherk ends that are not glued, two upwards and two
downwards.  More precisely: Given a configuration $(\ell, \psi)$ consisting of
two finite sequences of the same length.  We expect, for sufficiently small
$\varepsilon > 0$, a family $\sM_\varepsilon$ of minimal surfaces of genus
$n-1$ in $\T_\varepsilon^2 \times \R$ with four Scherk ends, where the flat
$2$-torus
\[
  \T_\varepsilon^2 = \R^2 / 2\pi \langle
  \bT_1 (1 + \Lambda_1 \tau_\tm), \bT_2 (1 + \Lambda_2 \tau_\tm)
  \rangle.
\]
They lift to DPMSs $\tilde\sM_\varepsilon$ in $\R^3$ such that
\begin{itemize}
	\item For each $1 \le k < n$, there exists $H_k(\varepsilon)$ that diverges to
		infinity while $\sM_\varepsilon - (0,0,H_k(\varepsilon))$ converges to a
		doubly periodic Scherk surface $\sS_k$ as $\varepsilon \to 0$.  Moreover,
		$\varepsilon^{-2} (H_{k+1} - H_k) \to \ell_k$ as $\varepsilon \to 0$.

	\item As $\varepsilon\to 0$, $\tilde\sM_\varepsilon$ converges, after a
		scaling by $\varepsilon^2$, to foliations by vertical planes of the direct
		product of $\R^2$ with segments of length $\ell_k$, $1 \le k < n$.

	\item As $\varepsilon \to 0$, the phase difference between $\sS_{k+1}$ and
		$\sS_k$ converges to $\psi_k$.
\end{itemize}

The technique for constructing TPMSs applies with little change.  But the
balance condition is even more strict: $\ell$ and $\psi$ can only be sequences
of length $\le 2$.

\begin{theorem}\label{thm:KMR}
	The family $\sM_\varepsilon$ of doubly periodic minimal surfaces described
	above exists if and only if $\ell$ and $\psi$ are sequences of length $1$ or
	$2$.  If they are of length $1$, then $\sM_\varepsilon$ is the doubly
	periodic Scherk surface.  If they are of length $2$, then $\sM_\varepsilon$
	is the Karcher--Meeks--Rosenberg examples obtained by gluing two doubly
	periodic Scherk surfaces.
\end{theorem}

\medskip

The paper is organized as follows: We will first review the construction and
geometry of the doubly periodic Scherk surfaces in
Section~\ref{sec:doublyscherk}.  Section~\ref{sec:construction} is the major
part of this paper, where we will glue Scherk surfaces into TPMSs.  Many steps
in the construction are similar to~\cite{chen2024} hence will not be detailed.
Proposition~\ref{prop:alpharho} and Section~\ref{sec:balance} are however key
and new ingredients of the construction and will be treated very carefully.
Finally, we prove Lemma~\ref{lem:triangle} and Theorems~\ref{thm:Meeks}
and~\ref{thm:KMR}, in this order, in Section~\ref{sec:main}.  This then
concludes the proof of the Main Theorem~\ref{thm:main}.  The construction for
DPMSs is very similar to that of TPMSs, hence will only be mentioned in
Section~\ref{sec:main}.

\subsection*{Acknowledgement}

The authors thank Matthias Weber, Martin Traizet, and Peter Connor for helpful
discussions.  The first author thanks Gerd Schroeder-Turk for hosting him at
Murdoch University and for working together on the TPMS of genus 4 with doubly
Scherk limit.  Some figures in this paper are from Weber's website
\url{https://minimalsurfaces.blog/}.

\section{Doubly periodic Scherk surfaces}
\label{sec:doublyscherk}

\begin{figure}[hbt]
	\includegraphics[width=0.95\textwidth]{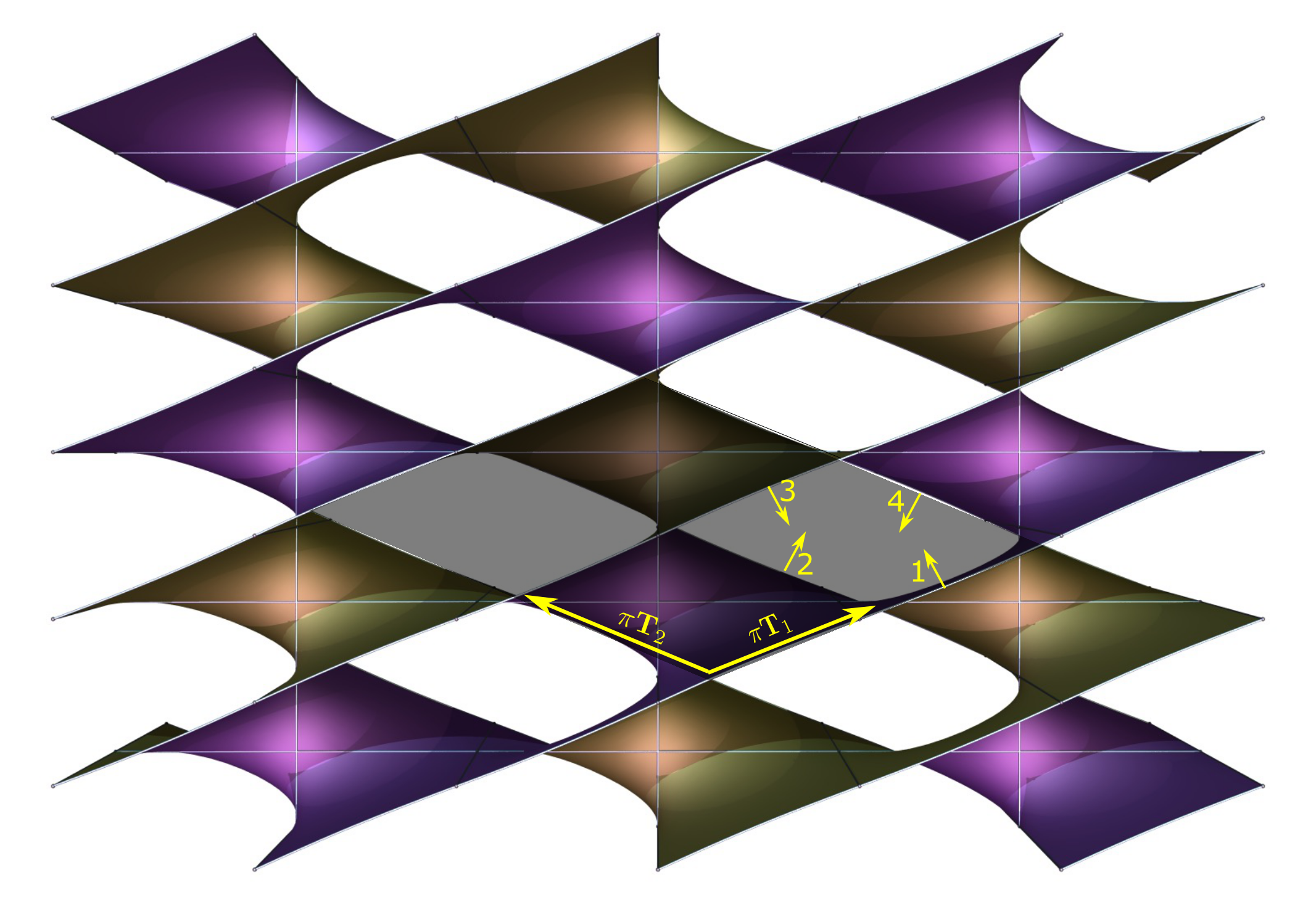}
	\caption{
    Top view of a doubly periodic Scherk surface.  The grey rhombus covers a
    fundamental domain of the translational symmetries, spanned by $2\pi\bT_1$
    and $2\pi\bT_2$.  The four ends and their normal vectors are marked with
    numbers, corresponding to the subscript $j$ in the text.
	}
	\label{fig:topview}
\end{figure}

A doubly periodic Scherk surface $\sS$ can be constructed as follows:  Let $P$
be a rhombus spanned by two vectors of length $\pi$, say,
\[
  \pi\bT_1=\pi(\cos\theta_1, \sin\theta_1)
  \qquad\text{and}\qquad
  \pi\bT_2=\pi(\cos\theta_2, \sin\theta_2).
\]
Up to a horizontal rotation, we may assume that $0 < \theta = \theta_1 = \pi -
\theta_2 < \pi/2$.  The Jenkins--Serrin Theorem~\cite{jenkins1966} guarantees
a minimal graph over $P$, unique up to vertical translations, that takes the
value $-\infty$ along the edges of $P$ parallel to $\bT_1$, and $+\infty$
along the edges parallel to $\bT_2$.  This minimal graph is bounded by the
four vertical lines over the vertices of $P$.  Order-$2$ rotations around
these vertical lines extend the graph into a DPMS, which is a \emph{doubly
periodic Scherk surface} $\sS$ with periods $2\pi\bT_1$ and $2\pi\bT_2$. See
Figure~\ref{fig:topview}.

\begin{figure}[hbt]
	\includegraphics[width=0.4\textwidth]{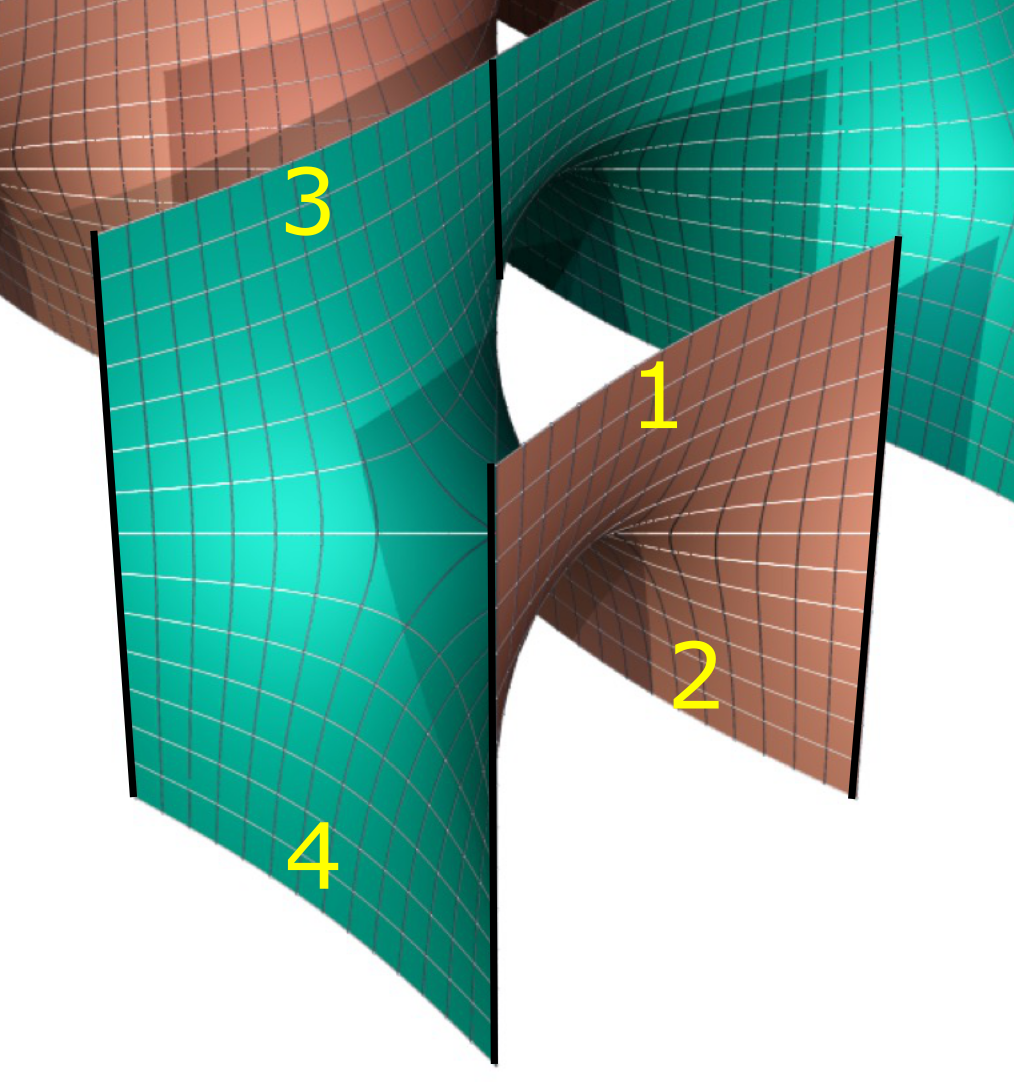}
	\caption{
    Part of a doubly periodic Scherk surface where the ends are marked and
    the vertical lines are highlighted.
	}
	\label{fig:standard}
\end{figure}

The Weierstrass parameterization of $\sS$ is given by
\[
  z \mapsto \int_0^z (\Phi_1, \Phi_2, \Phi_3),
\]
with the Weierstrass data
\[
  \Phi_1= \sum_{j=1}^4 \frac{-\ii\cos\theta_j}{z-p_j},\qquad
  \Phi_2= \sum_{j=1}^4 \frac{-\ii\sin\theta_j}{z-p_j},\qquad
  \Phi_3= \sum_{j=1}^4 \frac{\sigma_j}{z-p_j}dz
\]
defined on the Riemann sphere $\hat\C$ punctured at $p_j$, $1\le j \le 4$,
where $\sigma_j = -(-1)^j$, and $\theta_j$ are horizontal directions of the
ends with
\[
  \theta_1 = \theta,\quad
  \theta_2 = \pi - \theta,\quad
  \theta_3 = -\pi + \theta,\quad
  \theta_4 = -\theta,\quad
  0 < \theta < \pi/2.
\]
So the ends corresponding to $p_1$ and $p_3$ extend upwards and are parallel to
$\bT_1$, while the ends corresponding to $p_2$ and $p_4$ extend downwards and
are parallel to $\bT_2$.  See Figure~\ref{fig:standard}.

\begin{figure}[hbt]
	\includegraphics[width=0.4\textwidth]{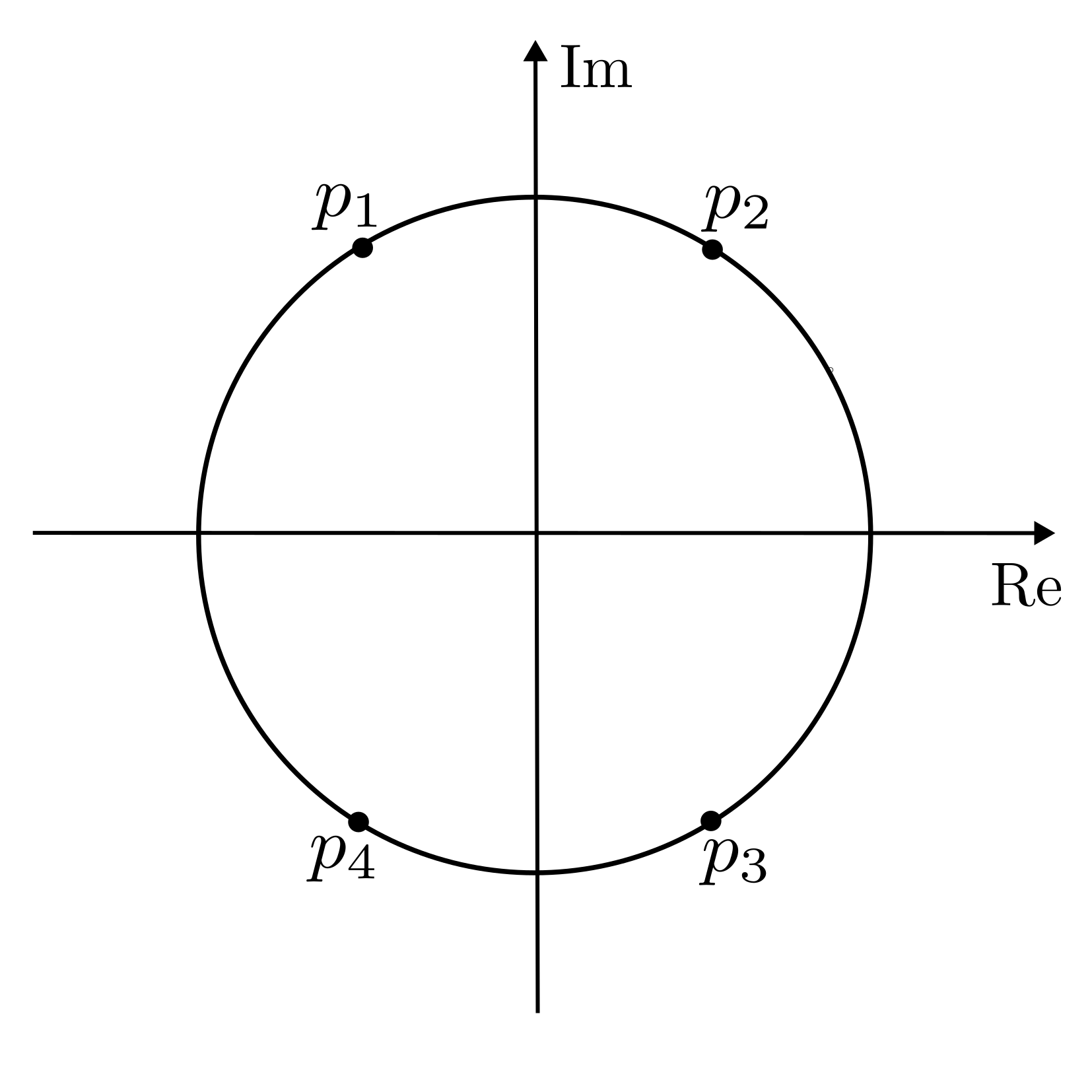}
	\caption{
    A doubly periodic Scherk surface is parameterized on a Riemann sphere with
    four punctures on the unit circle.
	}
	\label{fig:punctures}
\end{figure}

The order-$2$ rotations around vertical lines correspond to an anti-holomorphic
involution $\rho$ that preserves all the punctures.  For convenience, we assume
that the punctures all lie on the unit circle, then $\rho(z) = 1/\overline{z}$.
$\sS$ also contains two horizontal straight lines orthogonal to each other,
order-$2$ rotations around which interchange the upward and downward ends.
They can be assumed to correspond to reflections in the real and imaginary
axes.  So the punctures must be symmetrically placed on the unit circle.  See
Figure~\ref{fig:punctures}  We may choose
\[
  -1/p_1 = p_2 = 1/p_3 = -p_4 = \exp(\ii\vartheta).
\]
with $0 < \vartheta < \pi/2$.  Then the real line of $\hat\C$ is mapped to the
$y$-axis, and the imaginary line to the $x$-axis.  The conformality condition
$\Phi_1^2 + \Phi_2^2 + \Phi_3^3 = 0$ gives $\vartheta = \pi/2 - \theta$. Then
our choice leads to a very convenient (stereographically projected) Gauss map
\[
  G = -\frac{\Phi_1 + \ii \Phi_2}{\Phi_3} = z,
\]
which extends holomorphically to the punctures. So the normal vectors at the
ends are
\[
  N(p_j) = \sigma_j (-\sin \theta_j, \cos \theta_j, 0).
\]

Let $w_j = \ii (z-p_j)/(z+p_j)$ be a local coordinate around $p_j$.  Note that
$w_j$ is real positive on the circular segment between $p_j$ and $p_{j+1}$
(subscript modulo $4$) and
\[
  w_j \circ \rho(z)
  = \ii \frac{1/\overline{z} - p_j}{1/\overline{z} + p_j}
  = \ii \frac{\overline{p_j} - \overline{z}}{\overline{p_j} + \overline{z}}
  = \overline{w_j(z)}
\]
and, under this local coordinate, we have for $i=1$ and $2$,
\[
  \rho^* \Big( \frac{\Phi_i}{w_j} \Big) =
  - \frac{\overline\Phi_i}{\overline w_j},
\]
so $\Res{\Phi_i/w_j}{p_j} \in \ii\R$ is pure imaginary.

The four arcs between the punctures are mapped to the four vertical lines; see
Figure~\ref{fig:standard}.  Let
\[
  (\re\mu_j, \im\mu_j) :=
  \lim_{z \to p_j} \Big(
    -\cos(\theta_j) \arg(w_j) + \re\int_0^z \Phi_1,
    -\sin(\theta_j) \arg(w_j) + \re\int_0^z \Phi_2
  \Big)
\]
be the horizontal coordinates of the image of the arc between $p_j$ and
$p_{j+1}$ (subscript modulo $4$).  Then
\[
  \mu_1 = \pi\cos\theta,\quad
  \mu_2 = \pi\ii\sin\theta,\quad
  \mu_3 = -\pi\cos\theta,\quad
  \mu_4 = -\pi\ii\sin\theta,
\]
are at the vertices of a rhombus.

Expansion of the Weierstrass parameterization around $p_j$ in terms of $w_j$
gives
\begin{alignat*}{4}
	\re \int_0^z \Phi_1 &=
	\overbrace{\vphantom{\bigg(\bigg)}\re \mu_j + \cos\theta_j\arg(w_j(z))}^\text{planar terms}
	&&
	+ \overbrace{\re \left(w_j(z) \Res{\frac{\Phi_1}{w_j}}{p_j}\right)}^\text{undulation terms}
	&&+ O(|w_j(z)|^2),\\
	\re \int_0^z \Phi_2 &= \im \mu_j + \sin\theta_j\arg(w_j(z)) &&
	+ \re \left(w_j(z) \Res{\frac{\Phi_2}{w_j}}{p_j}\right)
	&&+ O(|w_j(z)|^2),\\
	\re \int_0^z \Phi_3 &= \phantom{\im} \nu_j+ \sigma_j \log|w_j(z)|&& &&
	+ O(|w_j(z)|),
\end{alignat*}
where
\[
  \nu_j = \lim_{z \to p_j}\Big( -\sigma_j \log|w_j(z)| + \re\int_0^z \Phi_3\Big)
  = -\sigma_j \log(\sin(\theta)\cos(\theta)).
\]

This shows that the ends are asymptotic to vertical planes, but also undulate
with an amplitude that decays exponentially as the height tends to $\pm\infty$.
See Figure~\ref{fig:standard}.  We then compute the ``initial'' amplitude in
the normal directions
\begin{align*}
	\Upsilon_j &= \sigma_j \Big\langle N(p_j), \Res{\frac{\Phi}{w_j}}{p_j}\Big\rangle_\mathrm{H}\\
	&= \ii \sin\theta_j \Res{\frac{\Phi_1}{w_j}}{p_j}
	-\ii \cos\theta_j \Res{\frac{\Phi_2}{w_j}}{p_j}
	=2.
\end{align*}
This could be more easily done using $G=z$ and noticing that
\[
  \Upsilon_j = \frac{\ii}{G}\frac{dG}{dw_j}(p_j).
\]

When we glue a sequence of Scherk surfaces, the $\sS$ constructed above will be
used to model odd-labeled Scherk surfaces $\sS_{2k-1}$.  As for the
even-labeled Scherk surfaces $\sS_{2k}$, we want that the ends corresponding to
$p_1$ and $p_3$ are still upwards but parallel to $\bT_2$, the ends
corresponding to $p_2$ and $p_4$ are still downwards but parallel to $\bT_1$,
and that the Gauss map coincides with $\sS$ at the punctures but has a pole at
$0$.  They will be modeled by a Scherk surface $\tilde\sS$ with the Weierstrass
data
\[
  \tilde\Phi_1= \sum_{j=1}^4 \frac{-\ii\cos\tilde\theta_j}{z-\tilde p_j},\qquad
  \tilde\Phi_2= \sum_{j=1}^4 \frac{-\ii\sin\tilde\theta_j}{z-\tilde p_j},\qquad
  \tilde\Phi_3= \sum_{j=1}^4 \frac{\tilde \sigma_j}{z-\tilde p_j}dz,
\]
where $\tilde\theta_j = \pi-\theta_j$, $\tilde\sigma_j = \sigma_j$, and $\tilde
p_j = p_j$.  Then the Gauss map $\tilde G = 1/G = 1/z$, and $G(p_j) = \tilde
G(\tilde p_{5-j})$ as we expect.  Hence the end of $\sS$ corresponding to $p_j$
will be glued to the end of $\tilde\sS$ corresponding to $\tilde p_{5-j}$.
Moreover, we will use the local coordinates $\tilde w_j = -\ii(z-\tilde
p_j)/(z+\tilde p_j)$.  This choice is consistent with the definition of phase.
Then we have
\[
  \tilde\Upsilon_j = \Upsilon_j = 2,\qquad
  \tilde\mu_j = \mu_{5-j},\qquad
  \tilde\nu_j = \nu_j.
\]
Note that here, $\tilde\mu_j$ gives the horizontal position of the arc between
$\tilde p_j$ and $\tilde p_{j-1}$.

\section{Construction of Triply Periodic Minimal Surfaces}\label{sec:construction}

In the following, all parameters vary around a central value.  The central
value of a parameter $x$ is denoted by an $\cent x$.  In particular, we fix a
configuration $(\cent\ell, \cent\psi)$ and three horizontal vectors
\[
  \bT_1 = \exp(\ii\theta_1),\quad
  \bT_2 = \exp(\ii\theta_2)
  \quad\text{and}\quad
  \bT_3 = \Psi_1\bT_1 + \Psi_2\bT_2,
\]
with $\theta_1 = \pi-\theta_2 = \theta$.  Moreover, we define
$\vartheta = \pi/2 - \theta$.

\subsection{Opening nodes}

Each Scherk surface $\sS_k$ is parameterized on a Riemann sphere $\hat\C_k$, $1
\le k \le n = 2m$, with four punctures on the unit circle corresponding to the
ends.  Unlike in the previous section, we use $p_{s, k}$ (resp.\ $q_{s, k}$),
$s \in \{+, -\}$, to denote the punctures corresponding to the upward (resp.\
downward) ends, and write $p = (p_{s, k})_{s=\pm, 1 \le k \le n}$ and $q=(q_{s,
k})_{s=\pm, 1 \le k \le n}$.

Our initial surface at $\varepsilon=0$ is a noded Riemann surface $\Sigma_0$
obtained by identifying punctures $p_{s,k}$ with $q_{s,k+1}$ for $1 \le k \le
n$, $s = \pm$.  Here and afterwards, the subscripts are taken modulo $n$.  The
central values for the punctures are defined by
\[
  -1/\cent p_{+,k} = \cent q_{-,k} = 1/\cent p_{-,k} = -\cent q_{+,k} = \exp(\ii\vartheta).
\]

As $\varepsilon$ increases, we open the nodes as follows: Depending on the
parity of $k$, consider local coordinates
\[
  u_{s,k} = (-1)^k \ii \frac{z-p_{s,k}}{z+p_{s,k}},\qquad
  v_{s,k} = (-1)^k \ii \frac{z-q_{s,k}}{z+q_{s,k}},
\]
in a neighborhood of the punctures. We fix a small number $\delta > 0$
independent of $k$ and $s$ such that, when the punctures are at their central
values, the disks $|\cent u_{s,k}| < 2\delta$ and $|\cent v_{s,k}| < 2\delta$ are
disjoint in each Riemann sphere $\hat\C_k$.  Then for punctures sufficiently
close to their central values, the disks $|u_{s,k}|<\delta$ and $|v_{s,k}|<\delta$
are disjoint in each $\hat\C_k$.

Consider complex parameters $t_{s,k}$, where $1 \le k \le n$ and $s = \pm$.  For
$t = (t_{s,k})$ in the neighborhood of $\cent t = 0$ with $0 < |t_{s,k}| <
\delta^2$.  We remove the disks
\[
  |u_{s,k}| < |t_{s,k}|/\delta, \qquad
  |v_{s,k+1}| < |t_{s,k}|/\delta
\]
and identify the annuli
\[
  |t_{s,k}| / \delta \le |u_{s,k}| \le \delta
  \quad \text{and} \quad
  |t_{s,k}| / \delta \le |v_{s,k+1}| \le \delta
\]
by
\[
  u_{s,k} v_{s,k+1} = t_{s,k},
\]
This produces a Riemann surface $\Sigma_t$.  Note that $\Sigma_t$ also depends
on the positions of the punctures, but the dependence is not written for
simplicity.

In the following, we consider the following fixed domains in all $\Sigma_t$:
\[
  U_{k,\delta}=\{ z \in \hat\C_k \colon |\cent u_{k,\pm}(z)| > \delta/2,|\cent v_{k,\pm}(z)| > \delta/2  \}
\]
and $U_{\delta}=\bigsqcup_{k = 1}^n U_{k,\delta}$. For $1 \le k \le n$, we
denote $0_k \in U_{k,\delta}$ the origin point in $\hat\C_k$, which will be
used as the starting point for the integration defining the Weierstrass
parameterization of the Scherk surface $\sS_k$. 

\subsection{Weierstrass data}

We construct a conformal minimal immersion using the Weierstrass
parameterization in the form
\[
  z\mapsto \re \int^z (\Phi_1, \Phi_2, \Phi_3) = \re \int^z \Phi ,
\]
where $\Phi = (\Phi_1, \Phi_2, \Phi_3) $ are holomorphic $1$-forms on
$\Sigma_t$.  At $t=0$, $\Phi$ extend holomorphically to so called \emph{regular
$1$-forms} on the noded Riemann surface $\Sigma_0$.  That means, they are
holomorphic away from the nodes, with simple poles at the punctures, and
opposite residues at $p_{s,k}$ and $q_{s,k+1}$.  Moreover, $\Phi$ are subject
to a conformal condition
\begin{equation}\label{eq:conformal}
	Q := \Phi_1^2 + \Phi_2^2 + \Phi_3^2 = 0.
\end{equation}
Note that $Q$ is a holomorphic quadratic differential on $\Sigma_t$.

\medskip

For the surface to be well defined, we need to solve the period problems.  We
will consider two groups of cycles on $\Sigma_t$.

For each $1 \le k \le n$ and $s = \pm$, let $A_{s,k}$ denote a small
counterclockwise circle in $U_{k,\delta}$ around $p_{s,k}$; it is homologous in
$\Sigma_t$ to a clockwise circle in $U_{k+1,\delta}$ around $q_{s,k+1}$.  Let
us write the periods in the vector form
\[
  \int_{A_{s,k}} \Phi = 2 \pi \ii (\ff_{s,k} - \ii \bx_{s,k}).
\]
By the Residue Theorem, it is necessary that 
\begin{align}
	\ff_{+,k} + \ff_{-,k} &= \ff_{+,k+1} + \ff_{-,k+1},\label{eq:balance}\\
	\bx_{+,k} + \bx_{-,k} &= \bx_{+,k+1} + \bx_{-,k+1}.\label{eq:xbalance}
\end{align}
By \cite[Proposition 4.2]{masur1976} (see also
\cite[Theorem~8.2]{traizet2013}), these equations guarantee the existence of
unique regular $1$-forms $\Phi(t)$ on $\Sigma_t$. Moreover, the restriction of
$\Phi(t)$ to $U_{\delta}$ depends holomorphically on $t$ in a neighborhood of
$0$.

The A-period problem requires that
\[
  \bx_{s,k} = s \bT_{\varsigma(k)} (1 + \Lambda_{\varsigma(k)} \tau_\tm)
\]
where $\bT_1$ and $\bT_2$ are the initial horizontal periods that we fixed at
the beginning of the section, and $\tau_\tm =
\exp(-\cent\ell_\tm/\varepsilon^2)$. So Eq.~\eqref{eq:xbalance} is
already solved.  Then $\Phi$ are uniquely determined by $\ff$ subject to
Eq.~\eqref{eq:balance}, which we refer to as the \emph{Balance Equation}.

More explicitly, at $\varepsilon=0$ and the central value of all parameters, we
have in $\hat\C_i$
\begin{align}
	\cent\Phi_1&=
	- \frac{\ii\cos(\theta_{\varsigma(k)})}{z-\cent p_{+,k}}dz
	- \frac{\ii\cos(\theta_{\varsigma(k+1)})}{z-\cent q_{-,k}}dz
	+ \frac{\ii\cos(\theta_{\varsigma(k)})}{z-\cent p_{-,k}}dz
	+ \frac{\ii\cos(\theta_{\varsigma(k+1)})}{z-\cent q_{+,k}}dz,\nonumber\\
	\cent\Phi_2&=
	- \frac{\ii\sin(\theta_{\varsigma(k)})}{z-\cent p_{+,k}}dz
	- \frac{\ii\sin(\theta_{\varsigma(k+1)})}{z-\cent q_{-,k}}dz
	+ \frac{\ii\sin(\theta_{\varsigma(k)})}{z-\cent p_{-,k}}dz
	+ \frac{\ii\sin(\theta_{\varsigma(k+1)})}{z-\cent q_{+,k}}dz,\nonumber\\
	\cent\Phi_3 &= 
	\frac{1}{z-\cent p_{+,k}}dz
	- \frac{1}{z-\cent q_{-,k}}dz
	+ \frac{1}{z-\cent p_{-,k}}dz
	- \frac{1}{z-\cent q_{+,k}}dz.\label{eq:explicitphi3}
\end{align}
In other words, $\cent\Phi$ is precisely the Weierstrass data of the Scherk
surface $\sS_k$ as we want.

\medskip

Let us now consider the period problems over the other group of cycles. For
every $1 \le k \le n$, $s = \pm$ and $t_{s,k}\neq 0$, let $B_{s,k}$ be the
concatenation of
\begin{enumerate}
	\item a path in $U_{k,\delta}$ from $0_k$ to $u_{s,k} = \delta$,

	\item the path parameterized by $u_{s,k} = \delta^{1-2\tau}\,(t_{s,k})^\tau$
		for $\tau \in [0,1]$, from $u_{s,k} = \delta$ to $u_{s,k} =
		t_{s,k}/\delta$, which is identified with $v_{s,k+1}=\delta$, and

	\item a path in $U_{k+1,\delta}$ from $v_{s,k+1}=\delta$ to $0_{k+1}$.
\end{enumerate}
Again for $1 \le k \le n$, let $B_k$ denote the concatenation
\[
  B_k = B_{+,k} * (-B_{-,k}),
\]
which is a cycle in $\Sigma_t$ that passes through two punctures.  We also need
to consider the cycle
\[
  B=B_{+,1} * B_{+,1} * \cdots * B_{+,n}.
\]
We need to solve the following B-period problem:
\begin{gather}
	\re \int_{B_k} \Phi = 0,\label{eq:Bk}\\
	\re \int_{B} \Phi = (\bT_3 + (\Psi_1 \Lambda_1 \bT_1 + \Psi_2 \Lambda_2 \bT_2) \tau_\tm, 1/\varepsilon^2)\label{eq:B},
\end{gather}
where $L$, $\bT_3$, $\Psi_1$ and $\Psi_2$ were defined in the Introduction.

\subsubsection{Dimension count}

Let us perform a dimension count before proceeding further. We have
\begin{itemize}
	\item $6n$ real parameters $\ff_{s,k}$.  

	\item $2n$ complex parameters $t_{s,k}$.  

	\item Up to M\"obius transforms, we may fix three punctures on $\hat\C_k$ for
		each $k$, leaving $n$ free complex parameters, say $p_{+,k}$.

	\item The shape parameters for the 3-torus.  Up to Euclidean similarities,
		the shape is determined by five parameters, namely $\varepsilon$,
		$\Lambda_1 - \Lambda_2$, $\Psi_1$, $\Psi_2$, and $\theta_1 - \theta_2$.
\end{itemize}
So we have $12n+5$ real parameters. Let us now count the equations.
\begin{itemize}
	\item The balancing equations~\eqref{eq:balance} consists of $3(n-1)$ real
		equations.

	\item The B-period problems~\eqref{eq:Bk} and~\eqref{eq:B} consist of $3(n+1)$
		real equations.

	\item The conformality equations~\eqref{eq:conformal} involve quadric
		differentials hence contain $3g - 3 = 3n$ complex equations.
\end{itemize}
Hence there are $12n$ real equations.  This is five less than the number of
parameters, so 5-parameter families are expected out of our construction.

\subsection{Using the Implicit Function Theorem}

\subsubsection{Conformality equations}

At $\varepsilon=0$ and the central value of all parameters, $\cent\Phi_3$ has
$4$ simple poles in $\hat\C_k$, hence $2$ zeros.  From the explicit
expression~\eqref{eq:explicitphi3}, we know that the two zeroes are
$\cent\zeta_{k,0} = 0_k$ and $\cent\zeta_{k,\infty} = \infty_k$.  When the
parameters are close to their central value, $\Phi_3$ has simple zeroes
$\zeta_{k,0}$ and $\zeta_{k,\infty}$ respectively close to $0_k$ and
$\infty_k$.

The conformality equation will be converted into the following equations:
\begin{align}
	\int_{A_{s,k}} \frac{Q}{\Phi_3} &= 0, \label{eq:conform1} \\
	\Res{\frac{Q}{\Phi_3}}{\zeta_{k,0}} &= 0. \label{eq:conform2}
\end{align}
for $1 \le k \le n$ and $s = \pm$.  The same argument as in~\cite{chen2024}
applies to prove the following statements.

\begin{proposition}\label{prop:alt-conformal}
	The conformality equation~\eqref{eq:conformal} is solved if the
	equations~\eqref{eq:balance}, \eqref{eq:conform1}, and \eqref{eq:conform2}
	are solved.
\end{proposition}

\begin{proposition}\label{prop:p}
	For $(t,\ff)$ in a neighborhood of $(0,\cent{\ff})$, there exists
	$p = (p_{s,k})$ depending analytically on $(t,\ff)$, such that the
	conformality equations~\eqref{eq:conform2} are solved, and
	$p_{s,k}(0,\cent{\ff}) = \cent p_{s,k}$.
\end{proposition}

From now on, we assume that the punctures $p$ are given by
Proposition~\ref{prop:p}.

For $j=1, 2$, let $\bT_j^\perp$ be the unit vector obtained by
anticlockwise rotation of $\bT_j$ by $\pi/2$.  We decompose
$\ff_{s,k}$ into
\[
  \ff_{s,k} = \alpha_{s,k} (\bT_{\varsigma(k)}, 0) + \beta_{s,k} (\bT_{\varsigma(k)}^\perp, 0) + \gamma_{s,k} (0, 0, 1).
\]
We also write $\rho_{s,k} = \sqrt{(\beta_{s,k})^2 + (\gamma_{s,k})^2}$ and
\begin{equation}\label{eq:t}
	t_{s,k} = -\exp\Big(-\ell_{s,k}\varepsilon^{-2} + \ii \psi_{s,k}\Big).
\end{equation}
The central values of $(\ell_{s,k}, \psi_{s,k})$ are $(\cent\ell_k,
\cent\psi_k)$.  We often change to the variable $\zeta_{s,k}$ where $\ell_{s,k} =
\cent\ell_k + \varepsilon^2\zeta_{s,k}$.  So $\cent\zeta_{s,k} = \partial \ell_{s,k}
/ \partial (\varepsilon^2)$ at $\varepsilon = 0$.

The next proposition is the key ingredient of our construction.  It is the only
place where a deformation in the horizontal lattice can play a role.  This is
the reason for our definition~\eqref{eq:horper} of the deformed horizontal
lattice.

\begin{proposition}\label{prop:alpharho}
	For $(t,\beta)$ in a neighborhood of $(0, 0)$, there exist unique values of
	$\alpha$ and $\rho$, depending real-analytically on $(t,\beta)$, such that
	the equations~\eqref{eq:conform1} are solved. At $t_{s,k}=0$ we have, no
	matter the values of other parameters, that
	\begin{equation}\label{eq:alphagamma1}
		\rho_{s,k} = 1
		\quad\text{and}\quad
		\alpha_{s,k} = 0.
	\end{equation}
	Moreover, at $(t,\beta)=(0,0)$, we have the Wirtinger derivatives
	\begin{equation}\label{eq:alphagamma2}
		\frac{\partial\rho_{s,k}}{\partial t_{s,k}} = \Xi_{s,k}
		\quad\text{and}\quad 
		\frac{\partial\alpha_{s,k}}{\partial t_{s,k}} = s \Xi_{s,k} \ii.
	\end{equation}
	where $\Xi_{s,k} = -2$ when $\cent\ell_k < \cent\ell_\tm$, and
	\[
    \Xi_{s,k} = -2-\frac{\Lambda_{\varsigma(k)}\exp(\cent\zeta_{s,k} - i\psi_{s,k})}{2}
	\]
	when $\cent\ell_{s,k} = \cent\ell_\tm$.
\end{proposition}

\begin{proof}
	Define for $1\le k \le n$
	\[
    \mathcal E_{s,k}(t,\alpha,\beta,\rho)=\frac{1}{2\pi\ii}\int_{A_{s,k}}\frac{Q}{\Phi_3}.
	\]
	Assume that $t_{s,k} = 0$. Then $\Phi_1$, $\Phi_2$, and $\Phi_3$ have a
	simple pole at $p_{s,k}$, so $Q/\Phi_3$ has a simple pole at $p_{s,k}$ and,
	by the Residue Theorem
	\[
    \mathcal E_{s,k}\mid_{t_{s,k}=0} =
    \frac{(\cent{\ff}_{k,s}-\ii \cent{\mathbf{x}}_{k,s})^2}{\gamma_{s,k}}=
    \frac{(\alpha_{s,k})^2+(\rho_{s,k})^2-2\ii s \alpha_{s,k}-1}{\gamma_{s,k}}.
	\]
	So the solution to~\eqref{eq:conform1} is $(\alpha_{s,k}, \rho_{s,k}) = (0, 1)$,
	no matter the values of the other parameters.  This
	proves~\eqref{eq:alphagamma1}.

	We compute the partial derivatives of $\mathcal E_{s,k}$ with respect to
	$\alpha_{s,k}$, $\rho_{s,k}$ at $(t,\alpha,\beta,\rho)=(0,0,0,1)$:
	\begin{equation} \label{eq:alphagamma4}
		\frac{\partial \mathcal E_{s,k}}{\partial \alpha_{s,k}} = -2\ii s,\qquad
		\frac{\partial\mathcal E_{s,k}}{\partial \rho_{s,k}} = 2.
	\end{equation}
	So the existence and uniqueness statement of the proposition follows from the
	Implicit Function Theorem.

	To prove the last point, we need to compute the partial derivative of
	$\mathcal E_{s,k}$ with respect to $t_{s,k}$ at $(t, \alpha_{s,k},
	\beta_{s,k}, \rho_{s,k}) = (0, 0, 0, 1)$.  The computation is similar to
	those in~\cite{chen2024}.  By the Residue Theorem
	and~\cite[Lemma~3]{traizet2008} (see also~\cite[Lemma~8.2]{chen2024}), we
	have
	\begin{align}
		\frac{\partial\mathcal E_{s,k}}{\partial t_{s,k}} =&
		\Res{\sum_{j=1}^3 2\frac{\cent\Phi_j}{\cent\Phi_3} \frac{\partial\Phi_j}{\partial t_{s,k}}}{\cent p_{s,k}}\nonumber\\
		=& -2\sum_{j=1}^3 \Res{\frac{\cent\Phi_j}{\cent\Phi_3}\frac{du_{s,k}}{(u_{s,k})^2}}{\cent p_{s,k}}
		\Res{\frac{\cent\Phi_j}{v_{s,k+1}}} {\cent q_{s,k+1}}\nonumber\\
		=& -2 \sum_{j=1}^3 \Res{ \frac{\cent\Phi_j}{u_{s,k}}}{\cent p_{s,k}} \Res{\frac{\cent\Phi_j}{v_{s,k+1}}}{\cent q_{s,k+1}} \label{eq:alphagamma5}\\
		&+ 2 \Res{\frac{\cent\Phi_3}{u_{s,k}}}{\cent p_{s,k}}
		\sum_{j=1}^3 \Res{\cent\Phi_j }{\cent p_{s,k}}\Res{\frac{\cent\Phi_j}{v_{s,k+1}}}{ \cent q_{s,k+1}}.\label{eq:alphagamma6}
	\end{align}
	Since $\cent Q=0$, we have
	\begin{align}
		\Res{\frac{u_{s,k} \cent Q}{d u_{s,k}}}{ \cent p_{s,k}}
		&= \sum_{j=1}^3 \Res{\cent\Phi_j}{ \cent p_{s,k}}^2 = 0, \label{eq:alphagamma7}\\
		\Res{\frac{\cent Q}{d v_{s,k}}}{ \cent q_{s,k}}
		&= 2\sum_{j=1}^3 \Res{\cent\Phi_j}{ \cent q_{s,k}} \Res{\frac{\cent\Phi_j}{v_{s,k}}}{\cent q_{s,k}}= 0. \label{eq:alphagamma8}
	\end{align}

	Because $\Res{\cent\Phi_j}{\cent q_{s,k+1}}= - \Res{\cent\Phi_j}{\cent
  p_{s,k}}$, Eq.~\eqref{eq:alphagamma8} implies that the
	term~\eqref{eq:alphagamma6} vanishes.  Note that
	\[
    \frac{1}{\sqrt{2}}\Res{\cent\Phi}{\cent p_{s,k}}, \quad
    \frac{1}{\sqrt{2}}\Res{\overline{\cent\Phi}}{\cent  p_{s,k}}, \quad
    \cent N( \cent p_{s,k})
	\]
	form an orthonormal basis of $\C^3$ for the standard hermitian product
	$\langle\cdot,\cdot\rangle_\text{H}$.  After recalling the definition and
	value of $\Upsilon$ in Section \ref{sec:doublyscherk}, we decompose in this
	basis
	\begin{align*}
		\Res{\frac{\cent\Phi}{u_{s,k}}}{\cent p_{s,k}} =&
		\frac{1}{2} \left\langle
		\Res{\cent\Phi}{\cent p_{s,k}},
		\Res{\frac{\cent\Phi}{u_{s,k}}}{\cent p_{s,k}}
		\right\rangle_{\!\text{H}}
		\Res{\cent\Phi}{\cent p_{s,k}}\\
		&+
		\left\langle
		\cent N(\cent p_{s,k}),
		\Res{\frac{\cent\Phi}{u_{s,k}}}{\cent p_{s,k}}
		\right\rangle_{\!\text{H}}
		\cent N(\cent p_{s,k})\\
		=&
		\frac{1}{2} \left\langle
		\Res{\cent\Phi}{\cent p_{s,k}},
		\Res{\frac{\cent\Phi}{u_{s,k}}}{\cent p_{s,k}}
		\right\rangle_{\!\text{H}}
		\Res{\cent\Phi}{\cent p_{s,k}} + 2 \cent N(\cent p_{s,k}).
	\end{align*}
	Here, the component on $\Res{\overline{\cent\Phi}}{\cent p_{s,k}}$ vanishes
	because of~\eqref{eq:alphagamma8}.  In the same way, after recalling that
	$\cent{N}(\cent p_{s,k})=\cent{N}(\cent q_{s,k+1})$,
	\begin{multline*}
		\Res{\frac{\cent\Phi}{v_{s,k+1}}}{\cent q_{s,k+1}} =\\
		\frac{1}{2} \left\langle \Res{\cent\Phi}{\cent p_{s,k}}, \Res{\frac{\cent\Phi}{v_{s,k+1}}}{\cent q_{s,k+1}} \right\rangle_{\!\text{H}} \Res{\cent\Phi}{\cent p_{s,k}}
		\;-\;2\cent N(\cent p_{s,k}).
	\end{multline*}

	Hence by Equation ~\eqref{eq:alphagamma5}
	\[
    \frac{\partial\mathcal E_{s,k}}{\partial t_{s,k}}=
    -2\left\langle
    \overline{\Res{\frac{\cent\Phi}{u_{s,k}}}{\cent p_{s,k}}},
    \Res{\frac{\cent\Phi}{v_{s,k+1}}}{\cent q_{s,k+1}}
    \right\rangle_{\!\text{H}}
    = 8.
	\]

	Since $\bx_{s,k}$ could also depend on $t_{s,k}$ (through $\tau_\tm$), we
	still need to compute
	\begin{equation}\label{eq:Ex}
		\Res{2\frac{\cent\Phi_1}{\cent\Phi_3} \frac{\partial\Phi_1}{\partial x_{s,k}}\frac{\partial x_{s,k}}{\partial t_{s,k}}}  {\cent p_{s,k}}
		+ \Res{2\frac{\cent\Phi_2}{\cent\Phi_3} \frac{\partial\Phi_2}{\partial y_{s,k}}\frac{\partial y_{s,k}}{\partial t_{s,k}}}  {\cent p_{s,k}},
	\end{equation}
	where $\bx_{s,k} = (x_{s,k}, y_{s,k}, 0)$.  The computation is dependent on
	$\cent\ell_k$.

	\medskip

	Recall that $\tau_\tm = \exp(-\cent\ell_\tm/\varepsilon^2)$.  When
	$\cent\ell_k < \cent\ell_\tm$, $\partial \bx_{s,k} / \partial t_{s,k} = 0$,
	hence~\eqref{eq:Ex} vanishes.  In this case, taking the total derivative of
	$\mathcal E_{s,k}$ against $t_{s,k}$ gives
	\begin{equation}\label{eq:alphagamma9}
		8-2is\frac{\partial \alpha_{s,k}}{\partial t_{s,k}}+2\frac{\partial \rho_{s,k}}{\partial t_{s,k}}=0.\\
	\end{equation}
	On the other hand, since $\mathcal E_{s,k}$ is holomorphic in $t_{s,k}$, we
	have
	\[
    2is\frac{\partial \alpha_{s,k}}{\partial t_{s,k}}+2\frac{\partial \rho_{s,k}}{\partial t_{s,k}}=0.\\
	\]
	Hence the Wirtinger derivatives
	\[
    \frac{\partial\rho_{s,k}}{\partial t_{s,k}} = -2
    \quad\text{and}\quad 
    \frac{\partial\alpha_{s,k}}{\partial t_{s,k}} = -2\ii s.
	\]

	\medskip

	When $\cent\ell_k = \cent\ell_\tm$, $\partial \bx_{s,k} / \partial \tau_\tm =
	s \Lambda_{\varsigma(k)} \bT_{\varsigma(k)}$ and
	\[
    \frac{\partial \tau_\tm}{\partial t_{s,k}} = \lim_{\varepsilon \to 0} \frac{\tau_\tm}{t_{s,k}} = -\exp(\cent\zeta_{s,k}-\ii\psi_{s,k}),
	\]
	So \eqref{eq:Ex} equals
	\begin{align*}
		& 2 \ii s \exp(\cent\zeta_{s,k}-\ii\psi_{s,k})\Lambda_{\varsigma(k)}\bT_{\varsigma(k)}\cdot
		\Bigg(
			\Res{\frac{\cent\Phi_1}{\cent\Phi_3} \frac{1}{u_{s,k}}}{\cent p_{s,k}},
			\Res{\frac{\cent\Phi_2}{\cent\Phi_3} \frac{1}{u_{s,k}}}{\cent p_{s,k}}
		\Bigg)\\
		=& 2 \ii s \exp(\cent\zeta_{s,k}-\ii\psi_{s,k})\Lambda_{\varsigma(k)}\bT_{\varsigma(k)}\cdot
		(- \ii s \bT_{\varsigma(k)})
		= 2 \exp(\cent\zeta_{s,k}-\ii\psi_{s,k}) \Lambda_{\varsigma(k)}.
	\end{align*}
	In this case, taking the total derivative of $\mathcal E_{s,k}$ against $t_{s,k}$
	gives
	\begin{equation}\label{eq:alphagamma10}
		8 + 2 \Lambda_{\varsigma(k)}\exp(\cent\zeta_{s,k} - i\psi_{s,k})
		- 2 i s \frac{\partial \alpha_{s,k}}{\partial t_{s,k}}
		+ 2 \frac{\partial \rho_{s,k}}{\partial t_{s,k}}
		=0.
	\end{equation}
	On the other hand, since $\mathcal E_{s,k}$ and $\tau_\tm$ are holomorphic in
	$t_{s,k}$, we have
	\begin{equation}\label{eq:alphagamma11}
		2 i s \frac{\partial \alpha_{s,k}}{\partial t_{s,k}}
		+ 2 \frac{\partial \rho_{s,k}}{\partial t_{s,k}}
		=0.
	\end{equation}
	Hence the Wirtinger derivatives
	\begin{align*}
		\frac{\partial \rho_{s,k}}{\partial t_{s,k}} &=
		-2-\frac{\Lambda_{\varsigma(k)}\exp(\cent\zeta_{s,k}-i\psi_{s,k})}{2} = \Xi_{s,k},\\
		\frac{\partial \alpha_{s,k}}{\partial t_{s,k}} &=
		s \Big(-2-\frac{\Lambda_{\varsigma(k)}\exp(\cent\zeta_{s,k}-i\psi_{s,k})}{2}\Big) \ii = s \Xi_{s,k} \ii.
	\end{align*}
\end{proof}
In other words, we have
\begin{align}
	\alpha_{s,k} &= 4 s \im t_{s,k} + o(|t_{s,k}|),\label{eq:alphaot}\\
	\rho_{s,k} &= \begin{cases}
		1 - 4 \re t_{s,k} + o(|t_{s,k}|), & \cent\ell_{s,k} < \cent\ell_\tm,\\
		1 - 4 \re t_{s,k} + \Lambda_{\varsigma(k)} \exp(\cent\zeta_{s,k}) |t_{s,k}| + o(|t_{s,k}|), & \cent\ell_{s,k} = \cent\ell_\tm. \label{eq:betaot}
	\end{cases}
\end{align}
Or more concisely
\[
  \alpha_{s,k} = -4s \sin(\psi_{s,k}) |t_{s,k}|, \quad
  \rho_{s,k} = 1 + K_{s,k} |t_{s,k}| + o(|t_{s,k}|),
\]
where
\[
  K_{s,k} = \begin{cases}
    4\cos(\psi_{s,k}) + \Lambda_{\varsigma(k)} \exp(\cent\zeta_{s,k}),& \cent\ell_k = \cent\ell_\tm,\\
    4\cos(\psi_{s,k}),& \cent\ell_k < \cent\ell_\tm,
  \end{cases}
\]

From now on, we assume that the parameters $\alpha$ and $\rho$ are given by
Proposition~\ref{prop:alpharho}.

\subsubsection{B-period problem}

\begin{proposition}\label{prop:minus}
	For $(\varepsilon, \ell_{+,k}, \psi_{+,k}, \beta_{+,k})$ in a neighborhood of
	$(0, \cent\ell_k, \cent\psi_k, 0)$, there exist unique values for
	$(\ell_{-,k})$, $(\psi_{-,k})$, and $(\beta_{-,k})$, depending smoothly on
	$(\varepsilon, \ell_{+,k}, \psi_{+,k})$, such that the B-period
	equations~\eqref{eq:Bk} are solved.  Moreover, at $\varepsilon = 0$, we have
	$\ell_{+,k} = \ell_{-,k}$, $\psi_{+,k} = -\psi_{-,k}$, and $\beta_{+,k} =
	\beta_{-,k}$.
\end{proposition}

\begin{proof}
	By~\cite[Lemma 1]{traizet2002} (see also Appendix of \cite{chen2024}), the
	difference
	\begin{equation}\label{eq:Bsk}
		\Big(\int_{B_{s,k}} \Phi \Big) - (\ff_{s,k} - \ii\bx_{s,k})\log t_{s,k} 
	\end{equation}
	extends holomorphically to $t_{s,k} = 0$. Moreover, its value at $t=0$ is
	equal to
	\begin{multline*}
		\lim_{z\to p_{s,k}}\bigg[\Big(\int_{0_k}^{z}\Phi\Big) - (\ff_{s,k} - \ii\bx_{s,k})\log u_{s,k}(z)\bigg]\\
		-\lim_{z\to q_{s,k+1}}\bigg[\Big(\int_{0_{k+1}}^{z}\Phi\Big) - (\ff_{s,k+1} - \ii\bx_{s,k+1})\log v_{s,k+1}(z)\bigg],
	\end{multline*}
	whose real part, by the computations in Section~\ref{sec:doublyscherk},
	equals
	\[
    (0, 0, -2\log(\sin\theta \cos\theta)).
	\]

	Recall that
	\[
    \ff_{s,k} =
    \alpha_{s,k} (\bT_{\varsigma(k)}, 0) +
    \beta_{s,k} (\bT_{\varsigma(k)}^\perp, 0) +
    \gamma_{s,k} (0, 0, 1)
	\]
	and that $\alpha_{s,k} = O(\exp(-\ell_{s,k}/\varepsilon^2))$.  Then, as $t \to 0$,
	we have
	\begin{align*}
		(\bT_{\varsigma(k)}, 0) \cdot (\int_{B_k} \Phi)
		&\to
		(\psi_{-,k} + \psi_{+,k}) +
		\lim_{\varepsilon \to 0}
		(\alpha_{-,k}\frac{\ell_{-,k}}{\varepsilon^2}-
		\alpha_{+,k}\frac{\ell_{+,k}}{\varepsilon^2})
		= (\psi_{-,k} + \psi_{+,k}),\\
		\varepsilon^2 (\bT_{\varsigma(k)}^\perp, 0) \cdot (\int_{B_k} \Phi)
		&\to
		(\beta_{-,k} \ell_{-,k} - \beta_{+,k} \ell_{+,k}),\\
		\varepsilon^2 (0, 0, 1) \cdot (\int_{B_k} \Phi) 
		&\to
		(\ell_{-,k} \gamma_{-,k} - \ell_{+,k} \gamma_{+,k}).
	\end{align*}
	The period problem requires all these to vanish at $\varepsilon = 0$.  This
	is solved at $\varepsilon=0$ with
	\[
    \ell_{-,k} = \ell_{+,k}, \qquad \beta_{-,k} = \beta_{+,k}, \qquad \psi_{-,k} = -\psi_{+,k},
	\]
	and the Proposition easily follows from the Implicit Function Theorem.
\end{proof}

From now on, we assume that the parameters $\ell_-$, $\psi_-$, and $\beta_-$ are
given by Proposition~\ref{prop:minus}.

\subsubsection{Balance condition}\label{sec:balance}

Recall that
\[
  K_k = \begin{cases}
    4\cos(\psi_k) + \Lambda_{\varsigma(k)} \exp(\cent\zeta_k),& \cent\ell_k = \cent\ell_\tm,\\
    4\cos(\psi_k),& \cent\ell_k < \cent\ell_\tm,
  \end{cases}
\]

\begin{proposition}\label{prop:phiell}
	For the balance equations~\eqref{eq:balance} to be solved at $\varepsilon =
	0$, it is necessary that the central values $\cent\ell_k$, $\cent\psi_k$ for
	$\ell_{+,k}$, $\psi_{+,k}$ only depend on the parity of $k$.  That is,
	\[
    \cent\ell_k=\cent\ell_{\varsigma(k)}, \qquad
    \cent\psi_k = \cent\psi_{\varsigma(k)}.
	\]
	If this is indeed the case, then for $\varepsilon$ in a neighborhood of $0$,
	there exist unique values for $\ell_{+,1}$, $\psi_+$ and $\beta_+$, depending
	smoothly on $\varepsilon$ and other parameters, such that
	\begin{equation}\label{eq:solB}
		\psi_{+,1} \to \cent\psi_1 = 2\Psi_1 / n,\quad
		\psi_{+,2} \to \cent\psi_2 = 2\Psi_2 / n,\quad
		\ell_{+,1} + \ell_{+,2} \to \cent\ell_1 + \cent\ell_2 = 2/n
	\end{equation}
	as $\varepsilon \to 0$, and the $B$-period problem~\eqref{eq:B} as well as the
	horizontal components of the balance equation~\eqref{eq:balance} are solved.
	Moreover, as $\varepsilon \to 0$, we have
	\begin{equation}\label{eq:zeta1}
		\zeta_{+,1} +\zeta_{+,2} = 4\log(\sin\theta \cos\theta).
	\end{equation}
\end{proposition}

\begin{proof}
	Define $\tau_k = \exp(-\cent\ell_k/\varepsilon^2)$.  We want $(\ell_+,
	\psi_+)(0) = (\cent\ell, \cent\psi)$.  From~\eqref{eq:alphaot},
	\eqref{eq:betaot}, Proposition~\ref{prop:minus} and the definition of $t$, we
	know that
	\begin{align*}
		\lim_{\tau_k \to 0} \frac{\alpha_{s,k}}{\tau_k} &= -4\sin(\psi_k)\exp(-\cent\zeta_{+,k}), \\
		\lim_{\tau_k \to 0} \frac{\rho_{s,k}-1}{\tau_k} &= K_k\exp(-\cent\zeta_{+,k}),
	\end{align*}
	both are independent of $s$.  Since the forces on the odd (resp.\ even)
	levels are decomposed into the same basis, the balance
	equations~\eqref{eq:balance} require, in particular, that
	\begin{equation}\label{eq:balanceparity}
		\frac{1}{\tau_k}(\alpha_{+,k}, \rho_{+,k} - 1) =
		\frac{1}{\tau_k}(\alpha_{+, \varsigma(k)}, \rho_{+, \varsigma(k)} - 1),
	\end{equation}
	These equations, depending on the parity of $k$, can be solved at
	$\varepsilon=0$ only if $\cent\ell_k = \cent\ell_{\varsigma(k)}$.  To see
	this, assume without loss of generality that $\cent\ell_k <
	\cent\ell_{\varsigma(k)}$, then $\tau_{\varsigma(k)} = o(\tau_k)$.
	Consequently, the left side of~\eqref{eq:balanceparity} extends to
	$\varepsilon=0$ with a nonzero value, while the right side tends to $0$.  Now, if $\cent\ell_k = \cent\ell_{\varsigma(k)}$, then~\eqref{eq:balanceparity} is
	solved at $\varepsilon = 0$ by $\cent\psi_k = \cent\psi_{\varsigma(k)}$.  We
	have thus proved the necessity that $\cent\ell_k$ and $\cent\psi_k$ both
	depend on the parity of $k$.

	Now assume that $\cent\ell$ and $\cent\psi$ both depend on the parity of $k$.
	By the same argument as in the previous proof, we have
	\[
    \Big( \int_B\Phi_1, \int_B\Phi_2, \varepsilon^2\int_B\Phi_3 \Big)
    \to \frac{n}{2} (\psi_{+,1} \bT_1, \psi_{+,2} \bT_2, \ell_{+,1} + \ell_{+,2}).
	\]
	So the B-period~\eqref{eq:B} is solved at $\varepsilon=0$ by the central
	values given in~\eqref{eq:solB}, and the proposition follows from the
	Implicit Function Theorem.  Moreover, as $\varepsilon \to 0$, we have
	\[
    \frac{\frac{2\varepsilon^2}{n}\int_B\Phi_3 - (\cent\ell_1 + \cent\ell_2)}{\varepsilon^2} \to\zeta_{+,1} +\zeta_{+,2} - 4\log(\sin\theta\cos\theta),
	\]
	which should vanish to close the periods.

	The balance equations~\eqref{eq:balance} require that the force $\ff$ is
	independent of $k$.  Consequently $\alpha$, $\beta$, and $\gamma$ all depend
	on the parity of $k$.  So the balance equations reduce to
	\begin{align}
		\alpha_{+,1} \bT_1 + \beta_{+,1} \bT_1^\perp
		&=
		\alpha_{+,2} \bT_2 + \beta_{+,2} \bT_2^\perp,\label{eq:balancehor}\\
		\gamma_{+,1} &= \gamma_{+,2}.\label{eq:balancever}
	\end{align}
	Note that $\bT_1$ and $\bT_2$ are linearly independent, and so are
	$\bT_1^\perp$ and $\bT_2^\perp$.  Hence for any $\alpha_{+,1}$ and
	$\alpha_{+,2}$, there exists unique $\beta_{+,1}$ and $\beta_{+,2}$ that
	solves the horizontal balance~\eqref{eq:balancehor}.  In fact, an elementary
	computation gives
	\[
    \beta_{+,1} = \frac{\alpha_{+,1}}{\tan(2\theta)} - \frac{\alpha_{+,2}}{\sin(2\theta)},
    \qquad
    \beta_{+,2} = \frac{\alpha_{+,1}}{\sin(2\theta)} - \frac{\alpha_{+,2}}{\tan(2\theta)}.
	\]
	Note that $\beta_{+,1}^2 - \beta_{+,2}^2 = \alpha_{+,2}^2 - \alpha_{+,1}^2$.
\end{proof}

From now on, we assume that the parameters $\ell_{+,1}$, $\psi_+$, and
$\beta_+$ are given by Proposition~\ref{prop:phiell}.

The remaining parameters are $\ell_{+,2}$ (or $\zeta_{+,2}$) and the lattice
shape parameters, including $\varepsilon$, $\theta_{1,2}$, $\psi_{1,2}$, and
$\Lambda_{1,2}$.  Usually, we would use the Implicit Function Theorem to solve
$\ell_{+,2}$ of $\zeta_{+,2}$ as a function of the shape parameters.  But such
a solution may not exist for a given $\Lambda_{1,2}$; see
Remark~\ref{rm:nosolution} later.  In the following, we will instead solve
$\Lambda_1 - \Lambda_2$ as a function of $\ell_{+,2}$ and other shape
parameters.

\begin{proposition}\label{prop:Lambda}
	If $\cent\ell_k$, $\cent\psi_k$ depend only on the parity of $k$ and $K_k \ne
	0$ for all $1 \le k \le n$, then $\cent\ell$ is necessarily constant for the
	balance equations~\eqref{eq:balance} to be solved.  If this is indeed the
	case, then for $\varepsilon$ in a neighborhood of $0$, there exists unique
	values for $\Lambda_2 - \Lambda_1$, depending smoothly on other parameters,
	such that the vertical component of the balance
	equation~\eqref{eq:balancever} is solved.
\end{proposition}

\begin{proof}
	Note that
	\[
    |\ff_{+,k}|^2 - 1 = \rho_{+,k}^2 + \alpha_{+,k}^2 - 1
    = 2 K_k \exp(-\zeta_{+,k}) \tau_k + o(\tau_k).
	\]
	The balance equation~\eqref{eq:balance} implies that
	\[
    \frac{|\ff_{+,1}|^2 - 1}{\tau_2}
    =
    \frac{|\ff_{+,2}|^2 - 1}{\tau_2}.
	\]
	But if $\cent\ell_\tm = \cent\ell_1 > \cent\ell_2$, $\tau_1 = o(\tau_2)$,
	then as $\varepsilon \to 0$, the left hand side vanishes but the right hand
	side is $K_2\exp(-\zeta_{+,2}) \ne 0$.  This proves that $\cent\ell$ must
	be constant.

	In this case, we now solve the vertical balance~\eqref{eq:balancever}. Note
	that
	\[
    \gamma_{+,k} - 1 = \sqrt{\rho_{+,k}^2 - \beta_{+,k}^2} - 1
    = K_k \exp(-\zeta_{+,k}) \tau_k + o(\tau_k). 
	\]
	Now that $\tau_1 = \tau_2 = \tau_k$, the vertical balance~\eqref{eq:balancever} is
	solved at $\varepsilon=0$ if
	\[
    \exp(-\zeta_{+,k}) K_k = 4 \cos(\cent\psi_k) \exp(-\zeta_{+,k}) + \Lambda_{\varsigma(k)}
	\]
	is a constant independent of $k$.  This implies that $\zeta_{+,k} =
	\zeta_{+,\varsigma(k)}$ depend on the parity of $k$ and
	\begin{equation}\label{eq:zeta2}
		4 \cos(\cent\psi_1) \exp(-\zeta_{+,1}) - \cos(\cent\psi_2) \exp(-\zeta_{+,2})
		+ \Lambda_1 - \Lambda_2 = 0.
	\end{equation}
	The proposition follows since this balance equation is linear in $\Lambda_2 -
	\Lambda_1$.
\end{proof}

\begin{remark}
	We are not able to conclude a constant $\cent\ell$ if $K_k=0$ for some $k$.
	We believe that in this case, the surfaces converge to KMR examples, so a
	future construction that glues KMR surfaces into TPMSs is expected.
\end{remark}

\begin{remark}\label{rm:nosolution}
	The equations~\eqref{eq:zeta1} and~\eqref{eq:zeta2} give a quadratic equation
	of $\exp(-\zeta_{+,2})$ at $\varepsilon=0$.  But the equation may not have a
	positive real solution.  This is why we are not able to solve $\zeta_{+,2}$
	as a function of the shape parameters.  For instance, for the Scherk limit oP
	surfaces, we have $\cent\psi_1 = 0$ and $\cent\psi_2 = \pi$ hence
	$\cos(\cent\psi_1) = -\cos(\cent\psi_2) = 1$; see Figure~\ref{fig:meeks}
	middle.   But then, if we insist $\Lambda_1 = \Lambda_2 = 0$ ,
	\eqref{eq:zeta2} would have no solution at all.  In other words, the shape of
	the lattice can not be arbitrary.  In the case of oP, the particular choice
	of phase forces a change in the aspect ratio of the horizontal lattice. 
\end{remark}

Finally, by the same argument as in~\cite{chen2024} and using
Eq.~\eqref{eq:Bsk}, one proves that the immersion given by the Weierstrass
Representation is a regular embedding and has the geometry described in the
Introduction.

\section{Proofs of main results}
\label{sec:main}

We now prove the Main Theorem~\ref{thm:main}.

\subsection{The special case of triangular lattice}

\begin{proof}[Proof of Lemma~\ref{lem:triangle}]
  Assume that the upward (resp.\ downward) ends of $\sS_k$ are parallel to
  unit vectors $\bT_k$ (resp.\ $\bT_{k-1}$), and that $\langle 2\pi \bT_{k-1},
  2\pi \bT_k \rangle$ are the periods of $\sS_k$.  Here we adopt the
  convention that $\bT_{k-1}$ and $\bT_k$ form an acute angle.  Since adjacent
  ends can not be parallel, we have $\bT_k \ne \bT_{k-1}$.  For the sake of
  contradiction, assume that the direction of ends does \emph{not} depend on
  the parity of $k$.  Then there must be three consecutive integers, say $i-1,
  i, i+1$, such that $\bT_{i-1}, \bT_i, \bT_{i+1}$ are pairwise non-parallel.

  Let $\T_0^2$ be the limit horizontal lattice as $\varepsilon \to 0$.
  Without loss of generality, we may assume that $\T_0^2 = \R^2 / \langle
  2\pi\bT_{i-1}, 2\pi\bT_i \rangle$.  But $2\pi\bT_{i+1}$ is also a period.
  That is, $\bT_{i+1}$ is an integer combination of $\bT_i$ and $\bT_{i-1}$.
  That is
  \begin{equation}\label{eq:integer}
    \bT_{i+1} = a \bT_{i-1} + b \bT_i, \qquad a, b \in \Z.
  \end{equation}
  Moreover, the horizontal fundamental parallelogram of $\sS_i$ and
  $\sS_{i+1}$ must have the same area.  So the angle $\angle\widehat{\bT_{i-1}
  \bT_i}$ must equal to the angle $\angle\widehat{\bT_i \bT_{i+1}}$.  Taking
  the cross product of~\eqref{eq:integer} with $\bT_i$, we obtain $a = -1$.
  But since $\bT_{i+1} \ne \pm\bT_{i-1}$ and they are all unit vectors, we
  have $0 < b < 2$, so $b=1$.

  We then conclude the relation $\bT_{i+1} + \bT_{i-1} = \bT_i$ among three
  unit period vectors.  This is only possible in a triangular lattice.
\end{proof}

So, if the limit horizontal lattice is not triangular, $\bT_k$ must depend on
the parity of $k$, which is the case considered in this paper.  

\begin{remark}\label{rem:triangle}
	In the special case of triangular limit horizontal lattice, the vertical
	planes can be parallel $e^{i \pi j/3}, j = 1, 2, 3$.  This leads to a wide
	variety of possibilities.  The simplest configuration is $\bT_k = e^{i \pi
	k/3}, k \in \Z$, which gives the doubly periodic Scherk limit of a recently
	discovered family of TPMSs of genus 4~\cite{markande2018}.  A similar
	argument as ours would prove that, whenever $\bT_k = \bT_{k'}$, we must have
	\[
		\cent\ell_k = \cent\ell_{k'} \quad \text{and} \quad
		\cent\psi_k = \cent\psi_{k'}.
	\]
	But with three possible directions, one may easily create very complicated
	sequences, even non-periodic infinite sequences.  More curious is an analogy
	of Proposition~\ref{prop:Lambda}.  More specifically, how would the surface
	behave if the horizontal lattice deforms away from triangular?  Although very
	interesting, we do not plan to discuss these examples in the current paper.
\end{remark}

In~\cite{chen2024}, the singly periodic Scherk limit of the surfaces
in~\cite{markande2018} also stood out as a curious special case, and the
behavior when the horizontal lattice deforms away from triangular was not
understood neither.

\subsection{Otherwise, no new examples}

\begin{proof}[Proof of Theorem~\ref{thm:Meeks}]
	We have shown that, for the construction to succeed, the periodic sequence of
	Scherk surfaces must actually have period $2$.  The situation is the same as
	if we wanted to construct TPMSs of genus 3 by gluing two Scherk surfaces.  We
	will prove that these TPMSs of genus 3 belong to Meeks' family, characterized
	by an orientation-reversing translation.

	The translation $\pi (\bx_{+,1} + \bx_{+,2}, 0)$ is orientation reversing for
	the Scherk surfaces.  Let $\varrho$ denote this transformation.  Note that
	$\varrho$ swaps $p_{\pm,k}$ and, respectively, $q_{\pm,k}$.  Up to M\"obius
	transforms, we may assume that $p_{+,k} = - p_{-,k}$.  Then $\varrho(z) =
	-1/\overline{z}$.  If our surfaces were assumed to be symmetric under
	$\varrho$, then we would automatically have
	\[
    \ff_{-,k} = \ff_{+,k}, \qquad t_{-,k} = \overline{t}_{+,k}.
	\]
	All other parameters can be found using the Implicit Function Theorem as we did
	in the previous section.

	But by the Implicit Function Theorem, the solution is unique.  Hence the
	solution we found in the previous section must coincide with the symmetric
	TPMSg3 solution.  In other words, the TPMSs we construct must be of genus $3$
	and have an orientation-reversing translational symmetry.  That is, they
	belong to the Meeks family.
\end{proof}

\begin{remark}
	We see from Remark~\ref{rm:nosolution} that, near the doubly periodic Scherk
	boundary, the 5-parameter Meeks manifold is not parameterized by the shape of
	the lattice.  More specifically, $\Lambda_{1,2}$ are replaced by $\partial
	\ell/ \partial(\varepsilon^2)$ as parameters of this manifold.
\end{remark}

Our construction for TPMSs can easily be adapted to construct DPMSs.  It
suffices to leave the four punctures $p_{\pm, n}$ and $q_{\pm, 1}$ ``free''.
That is, they are not identified with other punctures to form any node.
Equivalently, we can think of the situation as if $q_{\pm, 1}$ (resp.\ $p_{\pm,
n}$) were identified to some $p_{\pm, 0} \in \C_0$ (resp.\ $q_{\pm, n+1} \in
\C_{n+1}$) to form nodes but were never opened, i.e.
\[
	t_{s, 0} = t_{s, n} \equiv 0
\]
but $|t_{s, k}| > 0$ for $1 \le k < n$.  We are now ready to prove
Theorem~\ref{thm:KMR}.

\begin{proof}[Proof of Theorem~\ref{thm:KMR}]
	By the computation in the proof of~\ref{prop:alpharho}, since $t_{s, 0} =
	t_{s, n} \equiv 0$ we have
	\[
    \rho_{s, 0} = \rho_{s, n} \equiv 1, \qquad \alpha_{s, 0} = \alpha_{s, n} \equiv 0.
	\]
	By a similar argument as in the proof of~\ref{prop:phiell}, if $n > 2$, we
	must also have $t_{s, 2} = 0$, which is against our assumption.  This proves
	that $n \le 2$.  But in the cases $n=1$ (genus $0$) and $n=2$ (genus $1$),
	doubly periodic surfaces with Scherk ends have been fully
	classified~\cite{lazardholly2001, perez2005}, so the surfaces we construct
	must be a Scherk surface or a KMR example.
\end{proof}

\bibliography{References}
\bibliographystyle{plain}

\end{document}